\documentclass[fleqn]{article}
\usepackage{amsmath,amssymb,amsthm,xcolor,framed}
\usepackage[english]{babel}

\usepackage[title]{appendix}

\usepackage{graphicx}

\numberwithin{equation}{section}

\usepackage[bookmarks]{hyperref}

\newcommand{\abs}[1]{{\left\vert #1\right\vert}}
\newcommand{\R}{{\mathbb R}}

\newcommand{\e}{\varepsilon}

\newcommand{\loc}{{\rm loc}}

\def\Xint#1{\mathchoice
	{\XXint\displaystyle\textstyle{#1}}%
	{\XXint\textstyle\scriptstyle{#1}}%
	{\XXint\scriptstyle\scriptscriptstyle{#1}}%
	{\XXint\scriptscriptstyle\scriptscriptstyle{#1}}%
	\!\int}

\def\XXint#1#2#3{{\setbox0=\hbox{$#1{#2#3}{\int}$}
		\vcenter{\hbox{$#2#3$}}\kern-.5\wd0}}

\newtheorem{thm}{Theorem}[section]
\newtheorem{prop}[thm]{Proposition}
\newtheorem{lem}[thm]{Lemma}

\theoremstyle{definition}
\newtheorem{rem}[thm]{Remark}
\newtheorem*{rem*}{Remark}

\title{Interaction energies in nematic liquid crystal suspensions}
\date{\today}
\author{Lia Bronsard
	\thanks{Department of Mathematics and Statistics, McMaster University,
		{\footnotesize \href{mailto:bronsard@mcmaster.ca}{bronsard@mcmaster.ca}.}
	}
	\and 
	{Xavier Lamy
	\thanks{ Institut de Math\'ematiques de Toulouse, Universit\'e Paul Sabatier, Toulouse, France.
		{\footnotesize \href{mailto:xlamy@math.univ-toulouse.fr}{xlamy@math.univ-toulouse.fr}.}
	}}
\and 
{Dominik Stantejsky
	\thanks{Institut Élie Cartan de Lorraine, University of Lorraine.
		{\footnotesize \href{mailto:dominik.stantejsky@univ-lorraine.fr}{dominik.stantejsky@univ-lorraine.fr}.}
	}}
\and
{Raghavendra Venkatraman
	\thanks{Department of Mathematics, The University of Utah.
		{\footnotesize \href{mailto:raghav@math.utah.edu}{raghav@math.utah.edu}.}
	}}
}
\begin{document}

\maketitle

\begin{abstract}
We establish, as $\rho\to 0$, 
an asymptotic expansion
for the minimal Dirichlet 
energy of $\mathbb S^2$-valued
 maps
outside a finite number of 
%three-dimensional 
particles of size $\rho$ with 
fixed
centers $x_j\in\R^3$,
under general anchoring conditions at the particle boundaries.
Up to a scaling factor,
this expansion is of the form
\begin{align*}
E_\rho = \sum_j \mu_j -4\pi\rho \sum_{i\neq j} \frac{\langle v_i,v_j\rangle}{|x_i-x_j|} +o(\rho)\,,
\end{align*}
where $\mu_j$ is the minimal 
 energy after zooming in at scale $\rho$ around each particle, and $v_j\in\R^3$ is
a torque determined by the far-field behavior of the
corresponding single-particle minimizer.
The above expansion highlights
 Coulomb-like interactions 
between the particle
centers.
This
agrees with
 the 
 \textit{electrostatics analogy}
 commonly used in the physics literature for 
 colloid interactions in nematic liquid crystal.
That analogy was pioneered by Brochard and de Gennes in 1970,
 based on a 
 formal  
 linearization argument.
We obtain here for the first time
 a precise estimate
 of the energy error introduced 
 by this linearization procedure.
\end{abstract}

\section{Introduction}

We investigate
a mathematical model of interactions between colloid particles immersed in a nematic liquid crystal. 
Nematic liquid crystals are characterized by their orientational order:
one can think of elongated molecules which tend
to align along a common direction.
Each immersed particle distorts this alignment at long range, 
inducing interactions with the other particles.
When the sizes of the particles are
much smaller than the distances between them, 
the physics literature 
developes
an electrostatic analogy to describe their interactions, see \cite{brocharddegennes70,RNRP96,LPCS98} 
and the survey \cite[\S~2]{musevic19}.
That analogy relies on linearizing, away from the particles,
 the equations which describe nematic alignment at equilibrium.
Our main result gives an estimate of the error
 introduced by this linearization, 
 under precise modelling assumptions which we describe
next.

We use the simplest order parameter to describe the nematic phase: a unit vector $n\in\mathbb S^2$
indicating the direction of alignment.
A liquid crystal filling a domain $\Omega\subset\R^3$
 is 
  described by a map
$
n\colon\Omega\to\mathbb S^2$,
and 
we assume 
that its energy
 is 
given by
 \begin{align*}
E(n)=\int_{\Omega}|\nabla n|^2\, dx
+F(n_{\lfloor \partial\Omega})\,,
\end{align*}
for some $F\colon H^{1/2}(\partial\Omega;\mathbb S^2)\to [0,+\infty]$
which accounts for the anchoring of
 liquid crystal molecules at the domain boundary.
 Note that minimizing configurations satisfy the harmonic map equation $-\Delta n=|\nabla n|^2n$ in $\Omega$.

\begin{center}
\begin{figure}[h]
\centering
\includegraphics[scale=1.0]{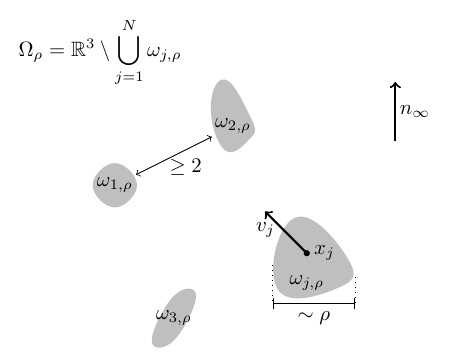}
\caption{General setup for Theorem~\ref{t:asympt}}
\label{fig:general_setting}
\end{figure}
\end{center}

Here we consider domains $\Omega$ and anchoring energies $F$ of a specific form, 
to model a system with $N$ foreign particles, 
all of the same small size $\rho>0$, 
but not necessarily the same shape, see Figure~\ref{fig:general_setting}. To be precise,
the
liquid crystal occupies the exterior domain
\begin{align*}
\Omega_\rho 
=\R^3\setminus \bigcup_{j=1}^N \omega_{j,\rho},
\qquad \omega_{j,\rho} =x_j +\rho\,\hat \omega_j\,,
\end{align*}
for fixed
  particle centers
   $x_1,\ldots,x_N\in\R^3$
and smooth open sets 
 %$\{\hat \omega_j\}_{j=1}^N$, each satisfying
\begin{align*}
\hat \omega_j\subset B_1 \subset\R^3\,,
\qquad
\text{for } j=1,\ldots,N\,.
\end{align*}
These open sets represent the particles after zooming in at scale $\rho$.

Rescaling by half the fixed minimal distance between these centers,
we assume without loss of generality that they
satisfy
\begin{align*}
|x_i-x_j|\geq 2\qquad\forall i\neq j \in\lbrace 1,\ldots,N\rbrace\,.
\end{align*}
We endow each
rescaled particle $\hat\omega_j$
 with an anchoring energy
\begin{align*}
&
\widehat F_j\colon H^{1/2}(\partial\hat\omega_j;\mathbb S^2)\to [0,\infty],
\quad\text{ weakly lower semicontinuous},
\end{align*} 
with non-empty domain 
$\lbrace \widehat F_j<\infty \rbrace\subset H^{1/2}(\partial\hat\omega_j;\mathbb S^2)$,
and  assume that anchoring at the boundary of each small particle
$\omega_{j,\rho}$ 
is described by
the rescaled energy
\begin{align*}
F_{j,\rho}(n_{\lfloor\partial\omega_{j,\rho}})=\widehat F_j(\hat n^\rho_{j\lfloor\partial\hat\omega_j}),\qquad
\hat n_j^\rho(\hat x) =n(x_j +\rho \, \hat x) 
\,.
\end{align*}
Examples of admissible anchoring energies $\widehat F_j$ are given in \cite[\S~1.2]{ABLV23}.
They include familiar examples of strong anchoring (Dirichlet conditions) 
and weak anchoring (enforced by a surface energy).
With these notations, the energy of a map 
$n\colon\Omega_\rho\to\mathbb S^2$ is given by
\begin{align}\label{eq:Erho}
E_\rho(n)=\frac{1}{\rho}\int_{\Omega_\rho}|\nabla n|^2\, dx 
+\sum_{j=1}^N 
F_{j,\rho}(n_{\lfloor\partial\omega_{j,\rho}})\,.
\end{align}
We impose  far-field alignment along a fixed 
orientation
 $n_\infty\in\mathbb S^2$ via the condition
\begin{align}\label{eq:farfieldOmegarho}
\int_{\Omega_\rho}\frac{|n(x)-n_\infty|^2}{1+|x|^2}\, dx <\infty\,.
\end{align}
Existence of a minimizer of $E_\rho$ under this far-field alignment constraint can be proved exactly as in \cite[\S~1.2]{ABLV23} for a single particle.
Our main result 
is an asymptotic expansion,
as $\rho\to 0$,
 of the minimal energy $E_\rho$.
That expansion depends on  minimizers of the single-particle problems
 \begin{align}\label{eq:muj}
&
\mu_j
=
\min
\bigg\lbrace
\widehat E_j(n)
\colon
\int_{\R^3\setminus\hat\omega_j}
\!\!
\frac{|n-n_\infty|^2}{1+|x|^2}\, dx <\infty
\bigg\rbrace\,,
\\
&\text{where }
\widehat E_j(n)
=
\int_{\R^3\setminus \hat\omega_j}|\nabla n|^2\, dx +\widehat F_j(n_{\lfloor\partial\hat\omega_j})
\,.
\nonumber
 \end{align}
 It is shown in \cite{ABLV23}
 that any minimizer 
 $\hat m_j$ of \eqref{eq:muj}
 has a far-field expansion
 \begin{align}\label{eq:vj}
 \hat m_j(x) =n_\infty +\frac{v_j}{|x|} +\mathcal O\Big(\frac{1}{|x|^2}\Big)\qquad
 \text{as }|x|\to +\infty\,,
 \end{align}
 for some $v_j\in \R^3$ orthogonal to $n_\infty$.
 The vector $v_j$ can be interpreted as a torque applied by the particle $\hat\omega_j$ on the nematic background \cite{brocharddegennes70}, see also \cite[Theorem~2]{ABLV23}.
The effective interaction between two particles depends on these vectors $v_j$.

\begin{thm}\label{t:asympt}
There exist  minimizers $\hat m_j$ of the single-particle problems
\eqref{eq:muj} such that the minimum of $E_\rho$ over maps $n\colon\Omega_\rho\to\mathbb{S}^{ 2}$ 
with far-field alignment \eqref{eq:farfieldOmegarho} satisfies
\begin{align}\label{eq:asympt}
\min
 E_\rho
&
=\sum_{j=1}^N \mu_j 
-4\pi\rho\sum_{i\neq j}
\frac{\langle v_i,v_j\rangle}{|x_i-x_j|}
+ o(\rho)\qquad\text{as }\rho\to 0\,,
\end{align}
where $\mu_j=\widehat E_j(\hat m_j)$ is the minimal single-particle energy \eqref{eq:muj},
and $v_j\in n_\infty^\perp$ is defined by the asymptotic expansion \eqref{eq:vj}
of $\hat m_j$.
\end{thm}

The interaction potential given by the term of order $\rho$ 
in the asymptotic expansion of Theorem~\ref{t:asympt}
corresponds to 
 solving
the Poisson equation with singular source term
\begin{equation} 
	\label{e.renormalizedenergy}
\Delta u_\rho = \sum_{j=1}^N 4\pi \rho v_j\delta_{x_j}\quad\text{ in }\R^3\,,
\quad
\text{that is,}
\quad
u_\rho(x) =\rho\sum_{j=1}^N \frac{v_j}{|x-x_j|}\,.
\end{equation}
The infinite energy of~$u_\rho$  
can indeed be renormalized to give
\begin{align*}
\lim_{\sigma\to 0}
\bigg(
\frac{1}{\rho}\int_{\R^3\setminus \bigcup B_\sigma(x_j)}\!|\nabla u_\rho|^2\, dx 
-\frac{4\pi}{\sigma}\rho\sum_{j=1}^N |v_j|^2\bigg)
=-4\pi\rho\sum_{i\neq j}\frac{\langle v_i,v_j\rangle}{|x_i-x_j|}
\,.
\end{align*}
This can be interpreted as follows:
\begin{itemize}
\item 
away from the particles,
the harmonic map equation $-\Delta n=|\nabla n|^2 n$ is linearized around the uniform state $n_\infty$,
which corresponds to writing
 $n\approx n_\infty +u$ and $-\Delta u\approx 0$;
\item
 the effect of the particle
 $\omega_{j,\rho}$ is replaced by a singular source term at $x_j$, and that source term is chosen to match the far-field expansion \eqref{eq:vj} generated by the single particle.
 \end{itemize}
This linearized description is the electrostatic analogy introduced in \cite{brocharddegennes70} and further developed in \cite{RNRP96,LPCS98}.
The difference in energy between the
(renormalized)
 linearized  description and the original nonlinear problem is what we estimate in
Theorem~\ref{t:asympt}.
This gives the asymptotic expansion \eqref{eq:asympt},
where
all the nonlinearity of the original problem 
is concentrated in
the presence of $\mu_j$ and $v_j$, 
determined by the single-particle problem \eqref{eq:muj}.

\subsubsection*{Ideas of proof} 

%As is standard for such results,
The proof of Theorem~\ref{t:asympt} consists of two  parts:
an upper bound, which we prove by constructing a competitor,
and a lower bound 
which we obtain via a precise description of minimizers $n_\rho$.

The competitor
we choose
 for the upper bound is 
equal to the single-particle minimizer~$\hat m_j$,
suitably rescaled,
in small regions $B_\sigma(x_j)$.
Outside these balls, we take its $\R^3$-valued harmonic extension 
	(tending to $n_\infty$ at far field) and project it back onto $\mathbb S^2$. 
	For a well-chosen $\sigma$ satisfying $\rho \ll\sigma\ll 1$, the energy of the competitor is controlled by the right-hand side of our expansion \eqref{eq:asympt}.

The lower bound is more challenging.
Thanks to classical compactness properties of energy-minimizing maps, 
  the blow-up at scale $\rho$ of a minimizer $n_\rho$ around each particle, given by
 $\hat n_j^\rho(\hat x)=n_\rho(x_j +\rho\,\hat x)$, converges  in 
 $H^1_{\loc}$ to a single-particle minimizer $\hat m_j(\hat x)$.
This provides the first term in the asymptotic expansion \eqref{eq:asympt}
and
%, in view of the structure of the upper bound construction,
%
%This 
suggests a natural route to obtain the next term:
  show that 
 $\hat n_j^\rho - \hat m_j$ is small enough in $B_{\sigma/\rho}$ to produce a negligible energy error, for an adequate scale $\sigma$.
If this were true, the conclusion would follow by
using  the energy of the harmonic extension of 
%the minimizer
 $n_\rho$ as a lower bound outside the regions $B_\sigma(x_j)$. 
In other words, that natural route
would require a quantitative rate for the convergence $\hat n_j^\rho\to\hat m_j$.
However, this convergence was obtained by weak compactness arguments, 
and quantifying it seems out of reach.
Instead, we modify our approach 
in order to conclude without quantitative rates.
This relies on the following two ingredients.
\begin{itemize}
\item 
The first  is a compensation effect between the inner and outer regions:
 if $\hat n_j^\rho$ is too different from $\hat m_j$ near $\partial B_{\sigma/\rho}$, the energy of the harmonic extension of $ n_\rho$ outside the regions $B_\sigma(x_j)$ 
 is increased
  by an amount which partly compensates the energy error inside $B_\sigma(x_j)$.
This 
 implies an improved lower bound for the full energy.
As a result, showing that the error is negligible
boils down to the estimate
$|\hat n_\rho^j(\hat x)-\hat m_j(\hat x)| \ll 1/|\hat x|$
 in the annulus $B_{2\sigma/\rho}\setminus B_{\sigma/\rho}$, 
for a choice of scale $\sigma \ll 1$ in an adequate range.
In terms of scaling, such estimate is consistent with the non-quantitative  $L^2$ convergence $\nabla\hat n_j^\rho\to\nabla\hat m_j$.
In comparison, 
separate lower bounds in the inner and outer regions,
 without taking advantage of this compensation effect, would have required
$|\hat n_\rho^j(x)-\hat m_j(x)| \ll \sqrt\rho/|\hat x|$ to make the error negligible.

\item 
The second ingredient
is a far-field expansion for $\hat n_j^\rho$ in large annuli, similar to the
expansion \eqref{eq:vj} of $\hat m_j$.
That far-field expansion eventually implies the estimate $|\hat n_\rho^j(x)-\hat m_j(x)| \ll 1/|\hat x|$,
 hence the conclusion thanks to the first ingredient.
The proof of the expansion \eqref{eq:vj} in \cite{ABLV23}
uses the fact that a classical harmonic function with finite energy in the exterior domain $\R^3\setminus B_\lambda$ only has radially decaying modes.
Here, in order to adapt it  to $\hat n_j^\rho$,
the main difference is that we must take into account
 radially increasing modes 
which can occur in an annulus,
and estimate them appropriately.
\end{itemize}

\subsubsection*{Related works}

Estimating the minimal energy of harmonic maps in exterior domains,
and interpreting it as an interaction energy,
is a very natural mathematical problem.
To the best of our knowledge, 
the perspective from which it has been addressed so far is different from the present one.
We wish to recall here 
the seminal works \cite{BCL} by Brezis, Coron and Lieb
in three dimensions,
and 
\cite[Chapter~I]{BBH} by Bethuel, Brezis and Hélein in two dimensions.
There, 
the objects of study are 
\textit{smooth}
$\mathbb S^2$ or $\mathbb S^1$-valued maps 
outside holes,
and the authors investigate the minimal energy within a fixed \textit{homotopy class}.
At first sight, 
their holes play a role very similar to our particles.
But here, 
on the contrary,
 our maps are not assumed to be smooth:
near the particles they may have several singularities, 
about which our analysis says nothing quantitative.
As a consequence, minimizing over a homotopy class would not even make sense in our setting,
and instead,
 admissible competitors are constrained by the anchoring conditions.
Finally, the results and methods in \cite{BCL} and \cite{BBH} 
are very different 
from each other 
but
 remain fundamentally nonlinear, 
while a linearization procedure is at the heart of the present work.

Note that in \cite{BBH}, the interaction energy is also obtained as the second term in an asymptotic expansion
%as the holes' size shrinks to zero, 
and is also of Coulomb type,
but this comes from the fact that $\mathbb S^1$-valued harmonic maps can be ``lifted'' to $\R$-valued harmonic maps, 
rather than a linearization around a uniform state as in the present work.
The analysis in
\cite{BBH} 
has initiated a  rich line of research, including generalizations to maps with values into general manifolds
and maps defined on higher-dimensional domains or manifolds, and we do not attempt here to give a list of these generalizations.

Finally, we mention the more recent papers~\cite{CS,GTVZ}. 
The paper~\cite{CS} uses methods from complex analysis and analogies with potential flows in fluid dynamics to study a version of our problem in the plane. 
The paper~\cite{GTVZ} considers interaction energies between particles in the so-called ``paranematic'' regime, in which nematic order is only felt at the boundaries of the particles. 
Consequently, the interaction energy is much more localized to essentially overlapping boundary layers, and the analysis is largely linear.

\subsubsection*{Further directions}

The physics of nematic suspensions 
raise many
mathematical questions,
and we mention here a few that
are directly linked with the present work.

We considered here the simplest model for the nematic phase.
Replacing the isotropic Dirichlet energy by a general anisotropic energy with three elastic constants \cite[\S~3.1.2]{degennes}
would likely be achievable at the cost of a few technical adjustments.
Adapting the present analysis to a 
$Q$-tensor model
(necessary to describe more symmetric single-particle minimizers, see e.g.  \cite{ABL16}), would require
new ingredients to deal with the extra length scale of phase transitions which is present in that model.

Recall that the vectors $v_j$ in \eqref{eq:asympt} 
can be interpreted as torques.
As detailed in \cite{ABLV23}, it follows from that interpretation that,
if the particles $\hat\omega_j$ are spherical, 
or if they are in an equilibrium orientation with respect to $n_\infty$, then all the vectors $v_j$ are zero.
In that case, our asymptotic expansion \eqref{eq:asympt} does not capture any interaction term.
These would be described by a next-order expansion, as predicted
by the electrostatics analogy \cite{LPCS98}.
An estimate of the error in that next order expansion would be very interesting.

From the physical point of view, it is also highly relevant to consider systems which are not at elastic equilibrium:
either because of other physical effects
(as already present in the original work of Brochard and de Gennes \cite{brocharddegennes70}
where the particles are magnetic),
or simply to describe time evolution.
The present work can serve as a first step towards 
these more complex models.

Finally, the limit $N\to \infty$ is of course very natural to study.
In that perspective, one goal could be to estimate the error in the continuum approximation proposed in \cite[\S~II.3]{brocharddegennes70}.
Another goal could be to establish a link between nematic suspensions and 
infinite point systems of Coulomb gas type, as has been done for Ginzburg-Landau vortices, see the survey \cite{serfaty14survey} and references therein.

\subsubsection*{Plan of the article}
In \S~\ref{s:harm_ext} we establish preliminary estimates on the energy of harmonic functions in exterior domains.
In \S~\ref{s:up} we 
present the upper bound construction.
In \S~\ref{s:low} we establish the lower bound, thus proving Theorem~\ref{t:asympt}.
In the Appendices \ref{a:poisson} and \ref{a:harm} 
we present for completeness some results about existence of decaying solutions to Poisson's equation, and estimates on the decay of harmonic functions in annuli.

\subsubsection*{Notations}
 
We write $A\lesssim B$ to denote $A\leq C B$ for a generic constant $C>0$, independent of $\rho$, but which can depend on the fixed parameters of our problem: $N$, $\hat\omega_j$, $\widehat F_j$, $n_\infty$.
We write $E_\rho(n;U)$ and $\widehat E_j(m,V)$ to denote restrictions of the integrals defining these energies to subsets $U\subset\Omega_\rho$ or $V\subset \R^3\setminus \hat \omega_j$. 
If  $\partial \omega_{j,\rho}\subset \partial U$ or $\partial\hat \omega_j\subset\partial V$, this notation is meant to include the corresponding anchoring term $F_{j,\rho}(n_{|\partial\omega_{j,\rho}})$ or $\widehat F_j(m_{|\partial\hat \omega_j})$.

\subsubsection*{Acknowledgements}
LB is supported by a NSERC Discovery grant.
XL is supported by the ANR project ANR-22-CE40-0006. 
A part of this work was carried out while DS was affiliated with McMaster University.
RV is supported by the U.S. National Science Foundation under the award NSF-DMS  2407592. 
A part of this work was conducted during a workshop at the Mathematisches Forschungsinstitut Oberwolfach in 2024.

\section{Prelude: harmonic extensions outside a union of small spheres}\label{s:harm_ext}

In this section we establish an estimate 
for the energy of harmonic functions $u$ in the
exterior domain 
\begin{align*}
U_\sigma=\R^3\setminus \bigcup_{j=1}^N B_\sigma(x_j),\qquad\text{for some }\sigma\in (\rho,1/2)\,,
\end{align*}
in terms of their boundary values on $\partial B_\sigma(x_j)$ that will be useful at several points in the proof of Theorem~\ref{t:asympt}.

We first introduce some notations.
We fix $\lbrace \Phi_k\rbrace_{k \in \mathbb{N}_0}$ an orthonormal Hilbert basis of $L^2(\mathbb S^2)$ which diagonalizes the Laplace-Beltrami operator,
\begin{equation*}
-\Delta_{\mathbb S^2}\Phi_k =\lambda_k \Phi_k,\qquad 0= \lambda_0 < \lambda_1\leq \lambda_2\leq \cdots \,. 
\end{equation*}
The set $\lbrace \lambda_k\rbrace_{k\in\mathbb N_0}$ coincides with $\lbrace \ell^2 + \ell\rbrace_{\ell\in\mathbb N_0}$, and the eigenfunctions corresponding to $\ell^2 + \ell$ span the homogeneous harmonic polynomials 
of degree $\ell$ (which have dimension $2\ell+1$). 
The function $\Phi_0$ is constant, equal to $1/(2\sqrt\pi)$.
Solutions $f(r)$ of $\Delta (f(r)\Phi_k(\omega))=0$ in $\R^3\setminus \lbrace 0\rbrace$  are spanned by $f_\pm(r)=r^{\pm\gamma_k^\pm}$, with
\begin{align*}
\gamma_k^+ &= \sqrt{\frac 14 +\lambda_k} - \frac{1}{2}   = \ell \quad \text{for }\lambda_k =\ell^2+\ell,\\
\gamma_k^- &= \sqrt{\frac 14 +\lambda_k} + \frac{1}{2}   = \ell +1 \quad \text{for }\lambda_k =\ell^2+\ell.
\end{align*}
These eigenfunctions satisfy the pointwise bound
\begin{align}\label{eq:ptwisePhik}
|\nabla^\alpha\Phi_k|
\lesssim 
\lambda_k^{\frac{1+\alpha}{2}}
\qquad\forall \alpha\geq 0\,,\: \: k\geq 1\,.
\end{align}
Indeed, for $k\geq 1$ and any $\omega_0\in\mathbb S^2$, 
in local coordinates around $\omega_0$ we can consider the rescaled function $\varphi (z) =\Phi_k(\omega_0 +z/\sqrt{\lambda_k})$ which satisfies 
$L \varphi  =\varphi $ in a fixed ball $B_1$ for some elliptic operator $L$ 
(with smooth coefficients depending on the local coordinates).
Elliptic estimates imply $|\varphi(0)|^2\lesssim \int_{B_1}|\varphi|^2\, dz \lesssim \lambda_k\int_{\mathbb S^2}|\Phi_k|^2$, hence $|\Phi_k|\lesssim \lambda_k^{1/2}$.
This shows \eqref{eq:ptwisePhik} for $\alpha =0$.
The case of higher derivatives $\alpha\geq 1$ follows from  
elliptic estimates on $\mathbb S^2$
and
the fact that  $(-\Delta)^\beta\Phi_k=\lambda_k^\beta\Phi_k$ for any integer $\beta\geq 1$.

\begin{prop}\label{p:harmextspheres}
If $\Delta u=0$ in $U_\sigma=\R^3\setminus \bigcup_{j=1}^N B_\sigma(x_j)$, $\int_{U_\sigma} |\nabla u|^2\, dx <\infty$ with $|u|\to 0$ as $|x|\to\infty$ and
$u=g_j$ on $\partial B_\sigma(x_j)$ for $j=1,\ldots ,N$,
with
\begin{align*}
g_j(x_j +\sigma\omega)=\sum_{k\geq 0}a_k^j \Phi_k(\omega),
\end{align*}
then
\begin{align}
\int_{U_\sigma}|\nabla u|^2\, dx
&
=
 \sigma\sum_j \sum_{k\geq 0} \gamma_k^-|a_k^j|^2
-\sigma^2\sum_j\sum_{i\neq j}\frac{ \langle a_0^i,a_0^j\rangle}{|x_i-x_j|}
\nonumber
\\
&\qquad
+\mathcal O(\sigma^3) \|a\|^2_{\ell^2}
\, .
\end{align}
\label{eq:harmextspheres}

\end{prop}

\begin{proof}[Proof of Proposition~\ref{p:harmextspheres}]

Consider $u_j$ the harmonic extension of $g_j(x_j +\cdot )$ in $\R^3\setminus B_\sigma$, given by
\begin{align}\label{eq:uj}
u_j(r\omega)
=\sum_{k\geq 0} a_k^j \left(\frac\sigma r\right)^{\gamma_k^-}\Phi_k(\omega),
\end{align}
and the function
\begin{align*}
\tilde u(x)=\sum_{j=1}^N u_j(x-x_j)\,,
\end{align*}
which is harmonic in $U_\sigma$.
The function $\tilde u - u$ is harmonic in $U_\sigma$ and satisfies $\tilde u - u=h_j$ on $\partial B_\sigma(x_j)$, where $h_j$ is smooth in $B_1(x_j)$ and given by
\begin{align}\label{eq:hj}
h_j(x_j+r\omega)
&
=\sum_{i\neq j} u_i(x_j-x_i +\sigma\omega)\,.
%\\
%&
%=h_j^0(x_j+r\omega) + h'_j(x_j+r\omega)\,.
%\nonumber
\end{align}
Since $u_i$ is decaying, 
this boundary error is  small for small $\sigma$,
and therefore,  its harmonic extension $u-\tilde u$ is also small.
Hence we expect that the energy of $u$ should coincide, at leading order, with the energy of $\tilde u$.
We will see that this heuristic is correct,
but we will also need to include next-order contributions 
to capture the second term in the right-hand side of \eqref{eq:harmextspheres}.
We 
start from the identity
\begin{align}\label{eq:errorutildeu}
\int_{U_\sigma}|\nabla u|^2\, dx
&
=
\int_{U_\sigma}|\nabla\tilde u|^2\, dx 
+\int_{U_\sigma} |\nabla\tilde u-\nabla u|^2\, dx
\nonumber
\\
&\quad
 +2\int_{U_\sigma}
\langle \nabla u -\nabla\tilde u,\nabla \tilde u\rangle \, dx\, .
\end{align}
The rest of the proof is structured as follows.
In Step~0 we
gather some  estimates on the boundary error $h_j$.
In Step~1 we estimate the last integral in \eqref{eq:errorutildeu}. The
 integral of $|\nabla\tilde u-\nabla u|^2$ is estimated in Step~2.
In Step 3 we compute $\int |\nabla\tilde u|^2$ and finally conclude in Step 4.

\medskip

\noindent\textit{Step~0. Estimates of the boundary error.}

\medskip

Let $i\neq j$, and $\alpha\geq 0$. 
Using that
$|x_i-x_j|\geq 2$,
$|\nabla_\omega^\alpha\Phi_k |
\lesssim \lambda_k^{(1+\alpha)/2}$
and $\gamma_k^-\lesssim 1+\sqrt{\lambda_k}$
we obtain, for all
 $x\in B_\sigma(x_j)$,
\begin{align*}
|\nabla^\alpha u_i(x - x_i)| 
&
\leq C_\alpha \sum_{k\geq 0}|a_k^i| 
\big(1+\lambda_k^{\frac{1+\alpha}{2}}\big) \sigma^{\gamma_k^-} 
\\
&
\leq C_\alpha \sigma \|a^i\|_{\ell^2} \bigg(\sum_{k\geq 0}
\big(1+\lambda_k^{1+\alpha}\big) \sigma^{2\gamma_k^--2}\bigg)^{1/2}\\
&
\leq C_\alpha \sigma \|a^i\|_{\ell^2}\,.
\end{align*}
The last inequality is valid because
 $\sigma \leq 1/2$ and $\gamma_k^-\geq 1$.
 In particular we have
\begin{align}
\label{eq:uixj}
\max_{0\leq\alpha\leq 4}
\sup_{B_\sigma(x_j)}
|\nabla^\alpha u_i(\cdot -x_i)| 
\lesssim
\sigma \|a^j\|_{\ell^2}\, .
\end{align}
Combining this bound with the definition \eqref{eq:hj} of the boundary error $h_j$, we infer
\begin{align}
|h_j(x_j+\sigma\omega)-h_j(x_j)|
&
\lesssim \sigma  \sum_{i\neq j} \sup_{B_\sigma(x_j)}|\nabla u_i(\cdot - x_i)|
\lesssim \sigma^2 
 \|a\|_{\ell^2}\,, \label{eq:hjdiff}
\end{align}
and
\begin{align}
|\Delta_\omega^2 \left[h_j(x_j+\sigma\omega)-h_j(x_j)\right]|
&
\lesssim \sigma^4 \sum_{i\neq j}
\max_{0\leq\alpha\leq 4} \sup_{B_\sigma(x_j)}|\nabla^4 u_i(\cdot - x_i)|
\nonumber
\\
&
\lesssim \sigma^5 
 \|a\|_{\ell^2} \label{eq:BiLaplhjdiff}
\, .
\end{align}

\medskip

\noindent\textit{Step 1. Estimating $\int\langle\nabla u-\nabla\tilde u,\nabla\tilde u\rangle$.}

\medskip

Since $\tilde u$ is harmonic in $U_\sigma$ and $\tilde u -u=h_j$ on $\partial B_\sigma(x_j)$,
we have
\begin{align*}
&
\int_{U_\sigma}
\langle \nabla u -\nabla\tilde u,\nabla \tilde u\rangle \, dx
=\int_{\partial U_\sigma} \langle u-\tilde u,\partial_\nu \tilde u\rangle \, d\mathcal H^2
\\
&
=\sigma^2\sum_j\int_{\mathbb S^2}\langle h_j(x_j +\sigma\omega),(\omega\cdot\!\nabla) \tilde u(x_j+\sigma\omega)\rangle\, d\mathcal H^2(\omega)
\\
&
=\sigma^2\sum_j 
\int_{\mathbb S^2}\langle h_j(x_j +\sigma\omega),\partial_r u_j(\sigma\omega)\rangle\, d\mathcal H^2(\omega)
\\
&
\quad
+\sigma^2\sum_j\sum_{\ell\neq j}
\int_{\mathbb S^2}
\langle h_j(x_j +\sigma\omega),(\omega\cdot\!\nabla)  u_\ell(x_j-x_\ell+\sigma\omega)\rangle\, d\mathcal H^2(\omega).
\end{align*}
To control the last integral we note that the estimate \eqref{eq:uixj} from Step~0 implies
\begin{align*}
&
|h_j(x_j+\sigma\omega)|
\lesssim \sigma N \|a\|_{\ell^2}
\quad
\text{and}
\quad
|\nabla u_\ell(x_j-x_\ell +\sigma\omega)|
\lesssim \sigma  \|a\|_{\ell^2}\,,
\end{align*}
for all $j\neq \ell$ and $\omega\in \mathbb S^2$.
We deduce
\begin{align}
&
\int_{U_\sigma}
\langle \nabla u -\nabla\tilde u,\nabla \tilde u\rangle \, dx
\nonumber
\\
&
=\sigma^2\sum_j 
\int_{\mathbb S^2}\langle h_j(x_j+\sigma\omega),\partial_r u_j(\sigma\omega)\rangle\, d\mathcal H^2(\omega)
%\\
%&\quad
+\mathcal O(\sigma^4)
\|a\|^2_{\ell^2}
\, .
\label{eq:errorutildeu2}
\end{align}
For all $j\in\lbrace 1,\ldots ,N\rbrace$, using the explicit expression \eqref{eq:uj} of $u_j(r\omega)$ we have
\begin{align}
&
\int_{\mathbb S^2}
\langle h_j(x_j+\sigma\omega),\partial_r u_j(\sigma\omega)\rangle\, d\mathcal H^2(\omega)
\nonumber
\\
&
=
-\frac{1}{\sigma}\sum_{k\geq 0}
\gamma_k^-
\int_{\mathbb S^2} \langle h_j(x_j+\sigma\omega),a_k^j \rangle \, \Phi_k(\omega)\, d\mathcal H^2(\omega)
\nonumber
\\
&
=
-\frac{1}{2\sqrt\pi\sigma}
\int_{\mathbb S^2} 
\langle h_j(x_j),a_0^j \rangle \, d\mathcal H^2
\nonumber
\\
&
\quad
-\frac{1}{\sigma}\sum_{k\geq 0}
\gamma_k^-
\big\langle a_k^j,
\int_{\mathbb S^2} 
\big(
h_j(x_j+\sigma\omega)-h_j(x_j)\big) \, \Phi_k(\omega)\, d\mathcal H^2(\omega) 
\big\rangle\, .
\label{eq:inthjdruj}
\end{align}
For the last 
equality we used the fact that the spherical harmonics $\Phi_k$ of order $k\geq 1$ have zero average, while $\Phi_0$ is constant equal to $1/(2\sqrt\pi)$.
To control
the last sum,
we note that by the estimate \eqref{eq:hjdiff} from Step~0 
 we have
\begin{align*}
\Big|
\int_{\mathbb S^2} 
\big(
h_j(x_j+\sigma\omega)-h_j(x_j)
\big) \, \Phi_0(\omega)\, d\mathcal H^2(\omega)
\Big|
&
\lesssim \sigma^2 
\|a\|_{\ell^2}\,,
\end{align*}
and, for $k\geq 1$, 
thanks to 
the fact that $\Phi_k=\lambda_k^{-2}\Delta_\omega^2\Phi_k$
and
the estimate
 \eqref{eq:BiLaplhjdiff} from Step~0,
\begin{align*}
&
\Big|
\int_{\mathbb S^2}
\big( 
h_j(x_j+\sigma\omega)-h_j(x_j)
\big) \, \Phi_k(\omega)\, d\mathcal H^2(\omega)
\Big|
\\
&
=\frac{1}{\lambda_k^2}
\Big|
\int_{\mathbb S^2}
\big( h_j(x_j+\sigma\omega)-h_j(x_j)\big) \, \Delta^2_\omega\Phi_k(\omega)\, d\mathcal H^2(\omega) 
\Big|
\\
&
=\frac{1}{\lambda_k^2}
\Big|
\int_{\mathbb S^2}
\Delta_\omega^2\big( h_j(x_j+\sigma\omega)-h_j(x_j)\big) \,  \Phi_k(\omega)\, d\mathcal H^2(\omega) 
\Big|
%\\
%&
 \lesssim \frac{\sigma^5 }{\lambda_k^2}
 \|a\|_{\ell^2}\,.
\end{align*}
From this and the previous inequality for $k=0$ we infer
\begin{align*}
&
\bigg|
\sum_{k\geq 0}
\gamma_k^-
\big\langle a_k^j,
\int_{\mathbb S^2} 
\big(
h_j(x_j+\sigma\omega)-h_j(x_j)\big) \, \Phi_k(\omega)\, d\mathcal H^2(\omega) 
\big\rangle
\bigg|
\\
&\qquad\qquad\qquad\lesssim 
\sigma^2 \|a\|_{\ell^2} |a_0^j|
+
 \sigma^5 \|a\|_{\ell^2}  \sum_{k\geq 1} \frac{\gamma_k^-}{\lambda_k^2}|a_k^j|
\lesssim 
\sigma^2\|a\|_{\ell^2}^2\,.
\end{align*}
The last inequality follows from the fact that $(\gamma_k^-/\lambda_k^2)^2\lesssim 1/\lambda_k^{3}$ is summable.
Using this to estimate the last line of \eqref{eq:inthjdruj} we deduce
\begin{align*}
\int_{\mathbb S^2}\langle h_j(x_j+\sigma\omega),&\partial_r u_j(\sigma\omega)\rangle\, d\mathcal H^2(\omega)
\\
&
=
-\frac{1}{2\sqrt\pi\sigma}
\int_{\mathbb S^2} 
\langle h_j(x_j),a_0^j \rangle \, \, d\mathcal H^2(\omega)
+\mathcal O(\sigma)
 \|a\|^2_{\ell^2}
\\
&
=
-\frac{ 2\sqrt\pi}{\sigma} \sum_{i\neq j}
u_i(x_j-x_i)
+\mathcal O(\sigma)
 \|a\|^2_{\ell^2}\,.
\end{align*}
The last equality
follows from the expression 
 \eqref{eq:hj} of $h_j$ and the fact that $|\mathbb S^2|=4\pi$.
Plugging this back into \eqref{eq:errorutildeu2} gives
\begin{align}\label{eq:prop21_mixed}
\int_{U_\sigma}
\langle \nabla u -\nabla\tilde u,\nabla \tilde u\rangle \, dx 
&
= 
-\sigma^2\sum_j \sum_{i\neq j}
\frac {2\sqrt\pi}\sigma \langle u_i(x_j-x_i),a_0^j\rangle
\nonumber
\\
&\quad
+\mathcal O(\sigma^3)
 \|a\|^2_{\ell^2}
\, .
\end{align}

\medskip

\noindent\textit{Step 2. Estimating $\int|\nabla \tilde u - \nabla u|^2$.}

\medskip

To bound the term $\int|\nabla \tilde u - \nabla u|^2$, we recall that $\tilde u-u$ is harmonic, apply Lemma~\ref{l:corrbdryharm} below and use \eqref{eq:uixj} to obtain
\begin{align}\label{eq:prop21_diff}
\int_{U_\sigma}|\nabla\tilde u -\nabla u|^2\, dx 
&
\lesssim \sigma \sum_j \|h_j\|^2_{C^1(\partial B_\sigma(x_j))}
\lesssim
\sigma^3 
\|a\|_{\ell^2}^2 \, .
\end{align}

\medskip

\noindent\textit{Step 3. Computing $\int|\nabla \tilde u|^2$.}

\medskip

Since $\tilde u$ is harmonic, we have
\begin{align*}
\int_{U_\sigma}|\nabla\tilde u|^2\, dx
&
=\int_{\partial U_\sigma} \langle \tilde u,\partial_\nu \tilde u\rangle \,d\mathcal H^2
=
\sigma^2\sum_{j,\ell,\ell'} I_j[u_\ell,u_{\ell'}],
\end{align*}
where, for $ j,\ell,\ell'\in\lbrace 1,\ldots,N\rbrace$,
\begin{align*}
I_j[u_\ell,u_{\ell'}]
&=
\frac{1}{\sigma^2} \int_{\partial B_\sigma(x_j)} 
 \langle u_\ell(\cdot -x_\ell),\partial_\nu  [u_{\ell'}(\cdot -x_{\ell'})]\rangle
 \,d\mathcal H^2
 \nonumber
\\
 &
=-
\int_{\mathbb S^2} \big\langle u_\ell(x_j-x_\ell +\sigma\omega),
(\omega\cdot\!\nabla)u_{\ell'}(x_j-x_\ell' +\sigma\omega)
\big\rangle
\, d\mathcal H^2(\omega).
\end{align*}
Since $u_\ell$ is small near $x_j-x_\ell$ for $j\neq \ell$,
the magnitude of this integral depends a lot on whether $\ell,\ell'$ are equal to $j$.
Main order terms will corresponding to 
 $\ell=\ell'=j$,
next-order to $\ell\neq \ell'=1$,
and all other terms will be negligible for our purposes.
We present next the estimates of each type of terms.

For $\ell=\ell'=j$, 
using the explicit expression \eqref{eq:uj} of $u_j(r\omega)$ 
we have
\begin{align*}
I_j[u_j,u_j]
&= -
\int_{\mathbb S^2}\langle u_j(\sigma\omega),\partial_r u_j(\sigma\omega)\rangle\, 
d\mathcal H^2(\omega)
%\\
%&
=
\frac 1\sigma\sum_{k\geq 0}\gamma_k^- |a_k^j|^2
\, .
\end{align*}
For $\ell\neq j$ and $\ell'=j$, using 
again  the explicit expression \eqref{eq:uj} of $u_j(r\omega)$,
and the fact that $\Phi_0=1/(2\sqrt\pi)$ while $\Phi_k$ has zero average for $k\geq 1$, 
we find
\begin{align*}
&I_j[u_\ell,u_j]
=
-\int_{\mathbb S^2} \big\langle u_\ell(x_j-x_\ell +\sigma\omega),\partial_r u_j(\sigma\omega)\rangle\, d\mathcal H^2(\omega)
\\
&
=\frac 1\sigma \sum_{k\geq 0}
\gamma_k^-
\int_{\mathbb S^2} \big\langle u_\ell(x_j-x_\ell +\sigma\omega),
a_k^j\rangle 
\,
\Phi_k(\omega)
\, d\mathcal H^2(\omega)
\\
&
=\frac{2\sqrt\pi}{\sigma}\langle u_\ell(x_j-x_\ell),a_0^j\rangle
\\
&
\quad
+\frac 1\sigma \sum_{k\geq 1}
\gamma_k^-
\int_{\mathbb S^2} \big\langle u_\ell(x_j-x_\ell +\sigma\omega)-u_\ell(x_j-x_\ell),
a_k^j\rangle 
\,
\Phi_k(\omega)
\, d\mathcal H^2(\omega)\,.
\end{align*}
The last line can be estimated as in Step~1 for the last sum in \eqref{eq:inthjdruj}, and we deduce
\begin{align*}
I_j[u_\ell,u_j]
&
=\frac{2\sqrt\pi}{\sigma}\langle u_\ell(x_j-x_\ell),a_0^j\rangle 
+\mathcal O(\sigma)\|a\|^2_{\ell^2}
\qquad\text{ for }\ell\neq j\,.
\end{align*}
For $\ell=j$ and $\ell'\neq j$ we find,
using 
  the explicit expression \eqref{eq:uj} of $u_j(r\omega)$,
and the estimate \eqref{eq:uixj} from Step~0,
\begin{align*}
I_j[u_j,u_{\ell'}]
&
=
\int_{\mathbb S^2} \big\langle u_j(\sigma\omega),(\omega\cdot\!\nabla) u_{\ell'}(x_j-x_{\ell'}+\sigma\omega)\big\rangle
\, d\mathcal H^2(\omega)
\\
&
=\sum_{k\geq 0}
\int_{\mathbb S^2} \big\langle a_k^j,(\omega\cdot\!\nabla) u_{\ell'}(x_j-x_{\ell'}+\sigma\omega)\big\rangle\,\Phi_k(\omega)
\, d\mathcal H^2(\omega)
\\
&
=\mathcal O(\sigma)\|a\|^2_{\ell^2}\, .
\end{align*}
Finally, for $\ell,\ell'\neq j$, we can directly use \eqref{eq:uixj} to deduce
\begin{align*}
I_j[u_\ell,u_\ell']
=\mathcal O(\sigma^2)\|a\|_{\ell^2}^2\,.
\end{align*}
Gathering all these estimates on the integrals $I_j[u_\ell,u_{\ell'}]$, we obtain
\begin{align}
\int_{U_\sigma}|\nabla \tilde u|^2\, dx
&
=\sigma^2\sum_{j,\ell,\ell'}I_j[u_\ell,u_{\ell'}]
\nonumber
\\
&
=
\sigma\sum_j \sum_{k\geq 0}\gamma_k^- |a_k^j|^2
+
 \sigma^2
 \sum_j \sum_{\ell\neq j}\frac{2\sqrt\pi}{\sigma} \langle u_\ell(x_j-x_\ell),a_0^j\rangle
 \nonumber
\\
&
\quad
+ \mathcal O(\sigma^3)
\|a\|_{\ell^2}^2\,.
\label{eq:prop21_tildeu}
\end{align}

\medskip

\noindent\textit{Step 4. Conclusion.}

\medskip

Inserting Equations \eqref{eq:prop21_mixed},\eqref{eq:prop21_diff} and \eqref{eq:prop21_tildeu} of Steps 1-3 into \eqref{eq:errorutildeu}, we end up with
\begin{align*}
\int_{U_\sigma}|\nabla\tilde u|^2\, dx
&
=
\sigma\sum_j \sum_{k\geq 0} \gamma_k^-|a_k^j|^2
-\sigma^2\sum_j\sum_{\ell\neq j}\frac {2\sqrt\pi}\sigma \langle u_\ell(x_j-x_\ell),a_0^j\rangle
\\
&\quad
+\mathcal O(\sigma^3) 
 \|a\|^2_{\ell^2}
\, .
\end{align*}
Finally, note that from \eqref{eq:uj} we find
\begin{align*}
\frac{2\sqrt\pi}{\sigma}u_\ell(x_j-x_\ell)=\frac{a_0^\ell}{|x_j-x_\ell|} +\mathcal O(\sigma)\|a\|_{\ell^2},
\end{align*}
which allows us to conclude.
\end{proof}

In Step~2 of Proposition~\ref{p:harmextspheres}'s proof
we used the following lemma to control the energy of a harmonic function with small boundary conditions.

\begin{lem}\label{l:corrbdryharm}
If $\Delta v=0$ in $U_\sigma$, $\int_{U_\sigma} |\nabla v|^2\, dx <\infty$ with $|u|\to 0$ as $|x|\to\infty$ and
$v=h_j$ on $\partial B_\sigma(x_j)$ for $j=1,\ldots ,N$,
with
\begin{align*}
h_j(x_j +\sigma\omega)=\sum_{k\geq 0}b_k^j \Phi_k(\omega),
\end{align*}
then, using the notation 
$\|b^j\|^2_{h^{1/2}} =\sum_{k\geq 0}(1+\sqrt{\lambda_k})|b_k^j|^2$, we have
\begin{align*}
\int_{U_\sigma}|\nabla v|^2\, dx\lesssim \sigma \sum_{j=1}^N \|b^j\|^2_{h^{1/2}}
\,.
\end{align*}
Moreover, if 
$h_j\in C^1(\partial B_\sigma(x_j))$
then
we have 
$\|b^j\|_{h^{1/2}}\lesssim \|h_j\|_{L^\infty} +\sqrt\sigma \|h_j\|_{C^1}$.
%
%$\|b^j\|_{h^{1/2}}\leq \|h_j\|_{L^\infty(\partial B_\sigma(x_j))} +\sqrt\sigma \|h_j\|_{C^1(\partial B_\sigma(x_j))}$.
\end{lem}

\begin{proof}
Denote by $v_j\colon\R^3\setminus B_\sigma$ the harmonic extension of $h_j(x_j +\cdot)$, that is,
\begin{align*}
v_j(r\omega)
=\sum_{k\geq 0} b_k^j \left(\frac\sigma r\right)^{\gamma_k^-}\Phi_k(\omega),
\end{align*}
so that, using the orthogonality of $\lbrace \Phi_k\rbrace$ and $\lbrace \nabla_\omega\Phi_k\rbrace$ in $L^2(\mathbb S^2)$, we have
\begin{align*}
\int_{\R^3\setminus B_\sigma}|\nabla v_j|^2\, dx
&
=
 \sum_{k\geq 0}((\gamma_k^-)^2+\lambda_k)|b_k^j|^2
\int_\sigma^\infty 
\left(\frac r{\sigma}\right)^{-2\gamma_k^-}\, dr
\\
&
=
\sigma \sum_{k\geq 0}\frac{(\gamma_k^-)^2+\lambda_k}{2\gamma_k^- -1}|b_k^j|^2
%\\
%&
\lesssim
\sigma \|b^j\|^2_{h^{1/2}},
\end{align*}
and
\begin{align*}
\int_{B_2\setminus B_1}|v_j|^2\, dx 
&
=\sum_{k\geq 0}\sigma^{2\gamma_k^-}|b_k^j|^2\int_{1}^2 r^{2-2\gamma_k^-}\, dr
\lesssim
\sigma^2 \|b^j\|_{\ell^2}.
\end{align*}
Next we fix a smooth cut-off function $\eta(r)$ such that $\mathbf 1_{r\leq 1}\leq \eta(r)\leq\mathbf 1_{r<2}$ and $|\eta'|\leq 2$ and set
\begin{align*}
\tilde v(x)=\sum_{j=1}^N \eta(|x-x_j|)v_j(x-x_j),
\end{align*}
so that $v=\tilde{v}$ on $\partial U_\sigma$ and by minimality of $v$ we have
\begin{align*}
\int_{U_\sigma}|\nabla v|^2\, dx 
&
\leq 
\int_{U_\sigma}|\nabla \tilde v|^2\, dx 
\leq N \sum_{j=1}^N \int_{B_2\setminus B_\sigma}|\nabla(\eta(r)v_j)|^2 \, dx
\\
&
\leq 2N \sum_{j=1}^N
\left(4\int_{B_2\setminus B_1}|v_j|^2\, dx
+ \int_{B_2\setminus B_\sigma}|\nabla v_j|^2
\right) \, dx\, .
\end{align*}
Combining this with the bounds on $v_j$ gives the estimate of $\int |\nabla v|^2\, dx $ in terms of $\|b_j\|_{h^{1/2}}$.

Assume moreover that $h_j\in C^1(\partial B_\sigma(x_j))$.
%From the identities
%\begin{align*}
%\int_{\partial B_\sigma(x_j)} |h_j|^2\,d\mathcal H^2
%&
%=
%\sigma^2 \sum_{k\geq 0} |b_k^j|^2,
%\\
%\int_{\partial B_\sigma(x_j)} |\nabla_\omega h_j|^2\, d\mathcal H^2
%&
%=
%\int_{\mathbb S^2} |\nabla_\omega h_j(x_j +\sigma\omega)|^2\, d\mathcal H^2(\omega)
%%\\
%%&
%=\sum_{k\geq 0}\lambda_k |b_k^j|^2,
%\end{align*}
For any $\alpha >0$ we estimate
\begin{align*}
2\sum_{k\geq 0}\sqrt{\lambda_k}|b_k^j|^2
\ &\leq \
\alpha \sum_{k\geq 0} |b_k^j|^2 +\frac 1\alpha \sum_{k\geq 0} \lambda_k |b_k^j|^2
\\
\ &\lesssim \
\frac{\alpha}{\sigma^2}\int_{\partial B_\sigma(x_j)} |h_j|^2 d\mathcal H^2
+
\frac1\alpha \int_{\partial B_\sigma(x_j)} |\nabla_\omega h_j|^2\, d\mathcal H^2
\\
\ &\lesssim \ \alpha \|h_j\|^2_{L^\infty(\partial B_\sigma(x_j))} +\frac{\sigma^2}{ \alpha}\|\nabla_\omega h_j\|^2_{L^\infty(\partial B_\sigma(x_j))}.
\end{align*}
We note that the claim of the lemma is trivial for constant $h_j$, so without loss of generality, we can assume that $\nabla_\omega h_j\not\equiv 0$, in particular $h_j\not\equiv 0$. 
This allows us to choose $\alpha=\sigma\|\nabla_\omega h_j\|_{L^\infty}/\|h_j\|_{L^\infty}$ gives
\begin{align*}
\sum_{k\geq 0}\sqrt{\lambda_k}|b_k^j|^2 
&
\lesssim
\sigma \|h_j\|_{L^\infty(\partial B_\sigma(x_j))}
\|\nabla_\omega h_j\|_{L^\infty(\partial B_\sigma(x_j))}
\\
&
 \lesssim 
\sigma \|h_j\|^2_{C^1(\partial B_\sigma(x_j))}.
\end{align*}
With $\|b^j\|_{\ell^2}\lesssim \|h_j\|_{L^\infty}$, this implies 
$\|b^j\|_{h^{1/2}}\lesssim  \|h_j\|_{L^\infty} +\sqrt\sigma \|h_j\|_{C^1}$.
\end{proof}

\section{Upper bound}\label{s:up}

In this section we perform the upper bound construction.

\begin{prop}\label{p:up}
The minimum of $E_\rho$ over all $n\colon\Omega\to\mathbb S^2$ 
with far-field alignment \eqref{eq:farfieldOmegarho} is bounded above by
\begin{align}\label{eq:up}
\min E_\rho \leq 
 \sum_j \mu_j
-\rho\sum_j\sum_{i\neq j} \frac{4\pi\langle v_i,v_j\rangle}{|x_i-x_j|}
+\mathcal O\left(\rho^{4/3}\right),
\end{align}
for any minimizers $\hat m_j$ of the
single-particle problem \eqref{eq:muj}, 
where
$\mu_j=\widehat E_j(\hat m_j)$
and $v_j$ 
is defined by the asymptotic expansion \eqref{eq:vj}.
\end{prop}

The upper bound is obtained by constructing a competitor and estimating its energy. 
In a ball $B_\sigma(x_j)$ around
each particle $\omega_{j,\rho}=x_j +\rho\widehat \omega_j$ we choose the competitor $n$
to be
equal to a single-particle minimizer $\hat m_j$, rescaled at scale $\rho$.
In the exterior $U_\sigma$ of these balls, 
we take $n$ to be the $\R^3$-valued harmonic extension, projected to $\mathbb S^2$.
The boundary values of this extension are determined by the maps $\hat m_j$, for which we have precise asymptotic estimates.
If $\sigma$ is large enough, the energy contribution inside each ball $B_\sigma(x_j)$ 
is close to $\mu_j=\widehat E_j(\hat m_j)$.
If $\sigma$ is not too large,
the energy contribution outside the balls $B_\sigma(x_j)$
can be accurately estimated using
Proposition~\ref{p:harmextspheres}.
Choosing $\sigma$ to balance error terms, we arrive at \eqref{eq:up}.

\begin{proof}[Proof of Proposition \ref{p:up}.]
We start by recalling from \cite{ABLV23} that for each $j=1,\ldots,N$, there exists a minimizer 
$ \hat m_j\colon\R^3\setminus\hat\omega_j\to\mathbb S^2$ 
 of $\widehat E_j$ under the far-field alignment constraint
 \begin{align*}
 \int_{\R^3\setminus\hat \omega_j}\frac{|\hat m_j-n_\infty|^2}{1+|x|^2}\, dx <\infty\,.
 \end{align*}
Furthermore, there exist $\lambda_0>0$
and $v_j\in n_\infty^\perp$ such that 
\begin{align}
&
\hat m_j(x)=
n_\infty +\frac{v_j}{r}
%-\frac{|v_j|^2}{2r^2}
%+ \sum_\alpha p_{j\alpha}\partial_\alpha \frac 1r 
+\hat w_j(x), \label{eq:mj-wj}
\\
&
|\hat w_j| +r|\nabla \hat w_j| +r^2 |\nabla^2 \hat w_j|
%+r^3 |\nabla^3 \hat w_j|
\lesssim \frac{1}{r^2}\qquad\text{for }r=|x|>\lambda_0 
\nonumber
\, .
\end{align}
Let $\sigma\in (\rho,1/2)$, to be fixed later on.
As explained above, we define our competitor to be equal to 
$\hat m_j((\cdot-x_j)/\rho)$ in each ball $B_\sigma (x_j)$.
At each boundary $\partial B_\sigma(x_j)$, it is therefore equal to  $n_\infty + g_j$,
where $g_j$ is given by
\begin{equation} \label{e.gj}
g_j(x_j +\sigma \omega)
= 
\hat m_j(\sigma \omega/\rho) -n_\infty 
= 
\frac{\rho}{\sigma} v_j +\hat w_j(\sigma\omega/\rho)
\, .
\end{equation}
We denote by $u$ the harmonic extension to $U_\sigma =\R^3\setminus \bigcup_j B_\sigma(x_j)$ satisfying $u=g_j$ on $\partial B_\sigma(x_j)$, as in Proposition~\ref{p:harmextspheres}.
With these notations, we define our competitor $n\colon\Omega_\rho\to\mathbb S^2$
by setting
\begin{align*}
n(x)=
\begin{cases}
\hat m_j\left(\frac{x-x_j}{\rho}\right)
&
\quad\text{if }|x-x_j|<\sigma,
\\
\frac{n_\infty +u}{|n_\infty +u|}
&\quad \text{if }x\in U_\sigma.
\end{cases}
\end{align*}

\begin{center}
\begin{figure}
\centering
\includegraphics[scale=1.0]{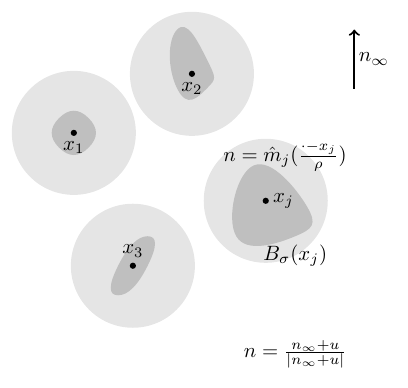}
\caption{Structure of the competitor $n$ constructed in Proposition~\ref{p:up}.}
\label{fig:upper_bound}
\end{figure}
\end{center}

We use the same notations as in Proposition~\ref{p:harmextspheres} and consider the spherical harmonics coefficients
\begin{align*}
a_k^j =\int_{\mathbb S^2}  g_j(x_j +\sigma\omega)\Phi_k(\omega)\,d\mathcal H^2(\omega)\,.
\end{align*}
Taking into account the decay properties \eqref{eq:mj-wj} of $\hat w_j$  we see that these coefficients $a_k^j$ satisfy
\begin{align*}
\Big|a_0^j -\frac{2\sqrt{\pi}\rho}{\sigma}v_j \Big|^2
\ &= \
\Big|
\frac {1}{2\sqrt\pi}
 \int_{\mathbb{S}^2} \Big(g_j(x_j+\sigma\omega) - \frac{\rho}{\sigma}v_j\Big)\,d\mathcal H^2(\omega) \Big|^2 \\
\ &= \
\Big|
\frac{1}{
2\sqrt{\pi}
}
\int_{\mathbb S^2} \hat w_j(\sigma\omega/\rho)\,d\mathcal H^2(\omega)
\Big|^2
\lesssim \frac{\rho^4}{\sigma^4}, \\
\sum_{k\geq 1} \lambda_k |a_k^j|^2 
\ &= 
\int_{\mathbb S^2} |\nabla_\omega [g_j(\sigma\omega)]|^2\, d\mathcal H^2(\omega)
\\
&
=
\int_{\mathbb S^2} |\nabla_\omega [\hat w_j(\sigma\omega/\rho)]|^2\, d\mathcal H^2(\omega)
%
%\
%\Big[g_j - \frac{\rho}{\sigma}v_j\Big]_{h^1}^2 
\\
&
 \leq \ 
\frac{\sigma^2}{\rho^2}\int_{\mathbb S^2}|\nabla \hat w_j|^2(\sigma\omega/\rho)\,d\mathcal H^2(\omega) 
\ \lesssim \
\frac{\rho^4}{\sigma^4} 
\, .
\end{align*} 
We deduce that for each~$j = 1,\ldots, N,$ 
\begin{align*}
a_j^0 &
=2\sqrt\pi\frac\rho\sigma v_j +\mathcal O(\rho^2/\sigma^2)\,,
\\
|a_j^0|^2 
&
= 4\pi \frac{\rho^2}{\sigma^2}|v_j|^2 +\mathcal O(\rho^3/\sigma^3)\,,
\\
\langle a_i^0,a_j^0\rangle
&
=4\pi\frac{\rho^2}{\sigma^2}\langle v_i,v_j\rangle
+\mathcal O(\rho^3/\sigma^3)\,,
\\
\sum_{k\geq 1}\gamma_k^- |a^j_k|^2
& 
\leq \sum_{k\geq 1}\lambda_k |a^j_k|^2 =\mathcal O(\rho^4/\sigma^4)
\,,
\\
\|a_j\|_{\ell^2}^2
&
=|a_j^0|^2 +\sum_{k\geq 1}|a_j^k|^2 
\leq |a_j^0|^2 +\sum_{k\geq 1}\lambda_k|a_j^k|^2
=\mathcal O(\rho^2/\sigma^2)\,.
\end{align*}
This enables us to estimate each term appearing in the asymptotic expansion \eqref{eq:harmextspheres} 
for the energy of $u$ provided by Proposition~\ref{p:harmextspheres}, namely
\begin{align*}
 \sigma\sum_j \sum_{k\geq 0} \gamma_k^-|a_k^j|^2
 &
 =
 \frac{\rho^2}{\sigma}\sum_j 4\pi|v_j|^2 + \mathcal O(\rho^3/\sigma^2)
 \\
 \sigma^2\sum_j\sum_{i\neq j}\frac{ \langle a_0^i,a_0^j\rangle}{|x_i-x_j|}
 &
 =
 \rho^2\sum_j 4\pi \frac{\langle v_i,v_j\rangle}{|x_i-x_j|}
 +\mathcal O(\rho^3/\sigma)
\\
\mathcal O(\sigma^3) \|a\|^2_{\ell^2}
&
=\mathcal O(\sigma\rho^2)\,.
\end{align*}
Dividing by $\rho$ and applying Proposition~\ref{p:harmextspheres}, we infer
\begin{align}
\frac 1\rho
\int_{U_\sigma}|\nabla u|^2\, dx
&
=
 \frac{\rho}{\sigma}\sum_j 4\pi|v_j|^2
-
\rho\sum_j 4\pi \frac{\langle v_i,v_j\rangle}{|x_i-x_j|}
\nonumber
\\
&\quad
 + \mathcal O(\rho^2/\sigma^2)
+\mathcal O(\sigma\rho)
\, .
\label{eq:energy_u_competitor}
\end{align}
Next we need to take into account the error introduced in that energy estimate by projecting $n_\infty +u$ onto $\mathbb S^2$.
To that end, note that 
the decay properties \eqref{eq:mj-wj} of $\hat w_j$
and the fact that $\langle n_\infty,v_j\rangle=0$
imply that the boundary condition 
$g_j(x_j+\sigma\omega)=(\rho/\sigma )v_j +\hat w_j(\sigma\omega/\rho)$ of $u$
at $\partial B_\sigma(x_j)$, 
defined in~\eqref{e.gj}, 
satisfies
$
\langle n_\infty,g_j\rangle =
\mathcal O(\rho^2/\sigma^2)
$.
So 
by the maximum principle
and since $u\to 0$ at $\infty$,
 we have
$\langle n_\infty ,u\rangle =\mathcal O(\rho^2/\sigma^2)$,
 and therefore
\begin{align*}
|n_\infty +u|^2 
&
\geq
1 +2\langle n_\infty,u\rangle 
 \geq  
1-C\frac{\rho^2}{\sigma^2}
\qquad\text{in }U_\sigma\,.
\end{align*}
Note also that for any smooth map $v$ with values in 
$\R^N\setminus\lbrace 0\rbrace$, 
the following inequality holds
\begin{align}
\Big|\partial_\alpha \Big[\frac{v}{|v|}\Big]\Big|^2
&
=\frac{1}{|v|^2}\Big|\partial_\alpha v -\big\langle \partial_\alpha v,\frac{v}{|v|}\rangle \frac{v}{|v|}\Big|^2
%\\
%&
=\frac{1}{|v|^2}
\Big(
|\partial_\alpha v|^2 - \big\langle \partial_\alpha v,\frac{v}{|v|}\big\rangle^2
\Big)
\nonumber
\\
&
\leq \frac{|\partial_\alpha v|^2}{|v|^2}\,.
\label{eq:norm_grad_proj}
\end{align}
Applying this to $v=n_\infty +u$,
we deduce
\begin{align*}
\frac 1\rho \int_{U_\sigma}|\nabla n|^2\, dx
 &\leq 
\frac{1}{\rho}\int_{U_\sigma} \frac{|\nabla u|^2}{|n_\infty +u|^2}\, dx
\\
 &\leq 
\frac{1}{\rho}\Big(1+C\frac{\rho^2}{\sigma^2}\Big)\int_{U_\sigma} |\nabla u|^2\, dx
\\
 &\leq 
\frac\rho\sigma \sum_j 4\pi|v_j|^2 
-\rho\sum_j\sum_{i\neq j} \frac{4\pi\langle v_i,v_j\rangle}{|x_i-x_j|}
+\mathcal O\left(\sigma\rho + \frac{\rho^2}{\sigma^2}\right)
\, .
\end{align*}
The last inequality comes from the estimate of $\int|\nabla u|^2$ in \eqref{eq:energy_u_competitor}.
In the whole domain, the energy of $n$ is therefore bounded by
\begin{align*}
E_\rho(n)
&
\leq \sum_j \widehat E_j(\hat m_j ; B_{\sigma/\rho}\setminus\hat \omega_j)
+
\frac\rho\sigma \sum_j 4\pi|v_j|^2 
\\
&
\quad
-\rho\sum_j\sum_{i\neq j} \frac{4\pi\langle v_i,v_j\rangle}{|x_i-x_j|}
+\mathcal O\left(\sigma\rho + \frac{\rho^2}{\sigma^2}\right)
\, .
\end{align*}
Noting that the expansion \eqref{eq:mj-wj} of $\hat m_j$ implies, for $\sigma\geq \lambda_0\rho$,
\begin{align}\label{eq:Emjextvj}
\widehat E_j(\hat m_j;\R^3\setminus B_{\sigma/\rho})
&
=\int_{|\hat x|\geq \frac\sigma\rho} |\nabla \hat m_j|^2\, d\hat x
=4\pi\frac\rho\sigma |v_j|^2+ \mathcal O\left(\frac{\rho^2}{\sigma^2}\right),
\end{align}
and recalling the definition $\mu_j
=
\widehat E_j(\hat m_j;\R^3\setminus \hat\omega_j)$, we deduce
\begin{align*}
\mu_j
 = 
\widehat E_j(\hat m_j;B_{\sigma/\rho}\setminus \hat\omega_j)
+ 4\pi\frac\rho\sigma |v_j|^2+ \mathcal O\left(\frac{\rho^2}{\sigma^2}\right)
\, .
\end{align*}
All together we obtain the upper bound
\begin{align*}
E_\rho(n)
&
\leq \sum_j \mu_j
-\rho\sum_j\sum_{i\neq j} \frac{4\pi\langle v_i,v_j\rangle}{|x_i-x_j|}
+\mathcal O\left(\sigma\rho + \frac{\rho^2}{\sigma^2}\right)
\, .
\end{align*}
Choosing
$\sigma=\sigma_\rho =\rho^{1/3}$
provides a 
remainder of order 
 $\rho^{4/3}$.
\end{proof}

\section{A matching lower bound}\label{s:low}

In this section we give the proof of a lower bound
which matches the upper bound of Proposition~\ref{p:up} at order $o(\rho)$ in the following sense:

\begin{prop}\label{p:low}
For any sequence $\rho\to 0$, there
exist minimizers $\hat m_j$ of the single-particle problem \eqref{eq:muj} and  a subsequence 
still denoted $\rho\to 0$ such that
\begin{align}\label{eq:low}
\min_{\eqref{eq:farfieldOmegarho}} E_\rho \geq 
 \sum_j \mu_j
-\rho\sum_j\sum_{i\neq j} \frac{4\pi\langle v_i,v_j\rangle}{|x_i-x_j|}
+o\left(\rho\right),
\end{align}
where $\mu_j=\widehat E_j(\hat m_j)$ is the minimal value of the single particle problem
\eqref{eq:muj} and $v_j$ is defined by the asymptotic expansion \eqref{eq:vj}.
\end{prop}

This proposition then allows us to prove Theorem~\ref{t:asympt}.
 
\begin{proof}[Proof of Theorem~\ref{t:asympt}]
Combining the upper bound of Proposition~\ref{p:up} 
 with the lower bound of Proposition~\ref{p:low}, we obtain the energy asymptotics
 \begin{align*}
& \min E_\rho =\sum_j \mu_j -4\pi \rho \, G +o(\rho),
\quad
\text{where }
G =\sup
\bigg\lbrace \sum_{i\neq j}
\frac{\langle v_i,v_j\rangle}{|x_i-x_j|}
\bigg\rbrace\,,
\end{align*}
and the supremum is over all collections of admissible vectors $v_j$, $j=1,\ldots,N$, which satisfy \eqref{eq:vj} for some minimizer $\hat m_j$.
For generic $n_\infty$, there is a unique admissible $v_j$  for each $j$, see \cite{ABLV23}, 
so the supremum is not needed
and this proves Theorem~\ref{t:asympt}.
If some $v_j$'s are not unique,
we need to check that the supremum in $G$ is attained to conclude
the proof of
 Theorem~\ref{t:asympt}.
To that end, it suffices to show that, for each particle $\hat \omega_j$, the set of admissible vectors $v_j$'s is compact.
This follows from two basic facts.
First, the set of minimizers $\hat m_j$ of the single-particle
problem \eqref{eq:muj} 
is compact in $H^1_{\loc}$ \cite{HKL88,luckhaus88}.
Second, the vector $v_j$ defined by \eqref{eq:vj} depends continuously on the minimizer $\hat m_j$ in that topology.
Assume indeed that $\hat m_j$ and $\tilde m_j$ are two minimizers with corresponding $v_j$ and $\tilde v_j$.
Then, using the asymptotic expansion \eqref{eq:mj-wj} 
for both minimizers 
we infer
\begin{align*}
|v_j-\tilde v_j|^2=R\int_{B_{2R}\setminus B_R}
|\nabla \hat m_j-\nabla \tilde m_j|^2\, dx +\mathcal O\Big(\frac{1}{R}\Big)\quad\text{as }R\to\infty.
\end{align*}
(It can be checked from the proof 
 \cite[Theorem~1]{ABLV23}
 of the expansion \eqref{eq:mj-wj} 
  that the constant in the estimate of the remainder can be taken independent from the minimizer.)
For any $\e>0$ we may choose $R$ large enough that the last term is smaller than $\e/2$, and we deduce that $|v_j-\tilde v_j|^2\leq \e$ provided $\hat m_j-\tilde m_j$ is small enough in
 $H^1(B_{2R}\setminus B_R)$, for this fixed radius $R$.
\end{proof}

In this whole section, we 
consider $n_\rho\colon\Omega_\rho\to\mathbb S^2$ a minimizer of $E_\rho$.
We extend $n_\rho$
 to $\mathbb R^3$ by filling the holes $\omega_{j,\rho}=x_j+\rho\hat\omega_j$ with $\mathbb S^2$-valued maps minimizing the Dirichlet energy,
and define the rescaled map
\begin{align}\label{eq:hatnjrho}
\hat n_j^\rho(\hat x)=n_\rho(x_j+\rho\hat x)\,,
\end{align}
around each particle $j\in\lbrace 1,\ldots,N\rbrace$.
We will freely 
extract subsequences and never make this explicit in the notations.

We divide the proof of Proposition~\ref{p:low} into four subsections.
In the first two we apply classical properties of energy-minimizing maps: small energy estimates and local $H^1$ compactness.
These provide, in
\S~\ref{ss:small_en_estim},  pointwise bounds that  will be used throughout the following sections;
and in \S~\ref{ss:comp_rescaled},
 strong $H^1_{\mathrm{loc}}$ convergence  
 of $\hat n_j^\rho$ to a minimizer $\hat m_j$ of the single-particle problem \eqref{eq:muj}.
Then, in \S~\ref{ss:lownjmj} 
we take advantage of a compensation effect to
 obtain, for any $1\ll\lambda\ll 1/\rho$,
 a lower bound 
which 
 depends
 on the smallness of $|\hat n_j^\rho -\hat m_j|$ 
on the annulus $B_{2\lambda}\setminus B_\lambda$.
Finally, in \S~\ref{ss:farfieldhatnrho} we establish a far-field expansion for $\hat n_j^\rho$ which 
we then combine with the far-field expansion \eqref{eq:vj} of $\hat m_j$ 
and 
the lower bound from the
 previous step 
 in order
 to
 conclude the proof of Proposition~\ref{p:low}.

\subsection{Pointwise estimates}\label{ss:small_en_estim}

In this section we gather pointwise estimates on 
$|\nabla\hat n_j^\rho|^2$ and $|\hat n_j^\rho -n_\infty|^2$
that follow from 
classical small energy estimates \cite{SU82} 
for the harmonic map $\hat n_j^\rho$ in the exterior of a large enough ball.

\begin{lem}\label{l:ptwise_estim}
There exists $\lambda_0>2$ such that,
 for all $\rho\in (0, 1/(2\lambda_0))$
and  $\lambda\in [\lambda_0,1/(2\rho)]$,
we have
\begin{align}
\label{eq:ptwise_estim_grad1}
\sup_{\partial B_\lambda} |\nabla \hat n_j^\rho|^2
&
\lesssim
\Xint{-}_{B_{5\lambda/4}\setminus B_{3\lambda/4}}
\!\!
|\nabla \hat n_j^\rho|^2 dx
\,,
\\
\lambda^2\sup_{\partial B_\lambda} |\nabla \hat n_j^\rho|^2
&
\lesssim
\Xint{-}_{B_{5\lambda/4}\setminus B_{3\lambda/4}}
\!\!
|\hat n_j^\rho-n_\infty|^2\, dx\,,
\label{eq:ptwise_estim_grad2}
\\
\label{eq:ptwise_estim_ninfty}
\sup_{\partial B_\lambda} |\hat n_j^\rho -n_\infty|^2
&
\lesssim 
\Xint{-}_{B_{5\lambda/4}\setminus B_{3\lambda/4}}
\!\!
|\hat n_j^\rho-n_\infty|^2\, dx\,.
\end{align}
\end{lem}
\begin{proof}[Proof of Lemma~\ref{l:ptwise_estim}]
Let $\eta>0$ be such that the small energy estimate
\begin{align}\label{eq:smallenSU}
&
\int_{B_1}|u-u_*|^2\, dx  + \int_{B_1}|\nabla u|^2 dx
 \leq \eta
 \nonumber
\\
&
\Rightarrow
\quad
\sup_{B_{1/2}}|\nabla u|^2  \lesssim 
\int_{B_1}|u-u_*|^2 dx\,,
\end{align}
is valid for any  map $u\colon B_1\to\mathbb S^2$ minimizing the Dirichlet energy 
with respect to its own boundary conditions, and any $u_*\in\R^3$  
(this is proved in \cite{SU82}, see \cite[\S~2.3]{simon96} for this precise statement).
Choosing $u_*$ 
in \eqref{eq:smallenSU} to be the average of $u$ on $B_1$, applying Poincaré's inequality, and 
decreasing $\eta$ if necessary, we also have the implication
\begin{align}\label{eq:smallenSU_avg}
&
\int_{B_1}|\nabla u|^2  dx
 \leq\eta
 \nonumber
 \\
 &
\Rightarrow
\quad
\sup_{B_{1/2}}|\nabla u|^2 \lesssim 
\int_{B_1}\Big|u-\Xint{-}_{B_1} u\Big|^2 dx
\lesssim \int_{B_1}|\nabla u|^2\, dx\,,
\end{align}
for any  map $u\colon B_1\to\mathbb S^2$ minimizing the Dirichlet energy.

\medskip

We introduce the notation $A_{\theta,\lambda}$ for the annulus of width $2\theta\lambda$ around the sphere $\partial B_\lambda$, 
that is,
\begin{align*}
A_{\theta,\lambda}=B_{(1+\theta)\lambda}\setminus B_{(1-\theta)\lambda},\qquad\text{ for }0<\theta<1\,.
\end{align*}
Let $\lambda_0>16$, to be chosen large enough later on.
We fix $\lambda\in [ \lambda_0, 1/(2\rho)]$ and $\theta\in [\lambda^{-1/4},1/2]$.
For any $x_0\in \partial B_\lambda$, we consider the  harmonic map
\begin{align*}
u(y)=\hat n_j^\rho(x_0+\theta\lambda y)\,,
\end{align*}
and check that it satisfies the smallness assumptions in \eqref{eq:smallenSU}-\eqref{eq:smallenSU_avg}.
We have
\begin{align*}
\int_{B_1}|\nabla u|^2\, dy 
&
=\frac {1}{\theta\lambda}
\int_{B_{\theta\lambda}(x_0)}|\nabla \hat n_j^\rho|^2\, dx
\leq 
\frac{1}{\theta\lambda} \int_{A_{\theta,\lambda}}|\nabla \hat n_j^{\rho}|^2\, dx
\\
&
\leq \frac{1}{\theta\lambda}E_\rho(n_\rho) \lesssim \frac{1}{\lambda^{3/4}}\,,
\end{align*}
since $\lambda^{1/4}\theta\geq 1$ and $E_\rho(n_\rho)\lesssim 1$,
so $u$ satisfies the smallness assumption in \eqref{eq:smallenSU_avg} provided $\lambda_0$ is large enough.
In addition, we have
\begin{align*}
\int_{B_1}|u-n_\infty|^2\, dy
&
=\frac{1}{(\theta\lambda)^3}\int_{B_{\theta\lambda}(x_0)}|\hat n_j^\rho-n_\infty|^2\, dx
\\
&
\leq
\frac{1}{\theta^3\lambda^3 }\int_{A_{\theta,\lambda}}
|\hat n_j^\rho -n_\infty|^2\, dx 
\\
&
\leq 
\frac{1}{\theta^3 \lambda}
\int_{\R^3\setminus B_1}
\frac{|\hat n_j^\rho -n_\infty|^2}{|x|^2}\, dx 
\lesssim 
\frac{1}{\lambda^{1/4}},
\end{align*}
since $\lambda^{3/4}\theta^3\geq 1$ and 
by Hardy's inequality combined with the bound $E_\rho(n_\rho)\lesssim 1$ we have
$\int |\hat n_j^\rho -n_\infty|^2/(1+|x|^2)\, dx\lesssim 1$, 
see also \eqref{eq:hardyhatnrhoj} later.
So $u$ satisfies the smallness assumption in \eqref{eq:smallenSU} with $u_*=n_\infty$, provided $\lambda_0$ is large enough.

\medskip

The estimate \eqref{eq:smallenSU_avg} thus gives
\begin{align*}
(\theta\lambda)^2 |\nabla \hat n_j^\rho(x_0)|^2 = |\nabla u(0)|^2\lesssim \frac{1}{\theta\lambda}
\int_{A_{\theta,\lambda}}|\nabla \hat n_j^{\rho}|^2\, dx\,,
\end{align*}
for any $x_0\in \partial B_\lambda$.
For $\theta=1/4$, this proves
\eqref{eq:ptwise_estim_grad1}.
The estimate \eqref{eq:smallenSU} with $u_*=n_\infty$ gives
\begin{align*}
(\theta\lambda)^2 |\nabla \hat n_j^\rho(x_0)|^2 = |\nabla u(0)|^2\lesssim 
\frac{1}{\theta^3\lambda^3 }\int_{A_{\theta,\lambda}}
|\hat n_j^\rho -n_\infty|^2\, dx\,,
\end{align*}
which, for $\theta=1/4$, proves \eqref{eq:ptwise_estim_grad2}.

\medskip

It remains to prove the pointwise estimate \eqref{eq:ptwise_estim_ninfty} on $|\hat n_\rho^j-n_\infty|$.
To that end we use the pointwise estimate \eqref{eq:ptwise_estim_grad2}
and the fundamental
 theorem of calculus to bound the oscillation of $\hat n_j^\rho$ on
the annulus $A_{\theta/2,\lambda}$, whose diameter is $\leq 4\lambda$. 
Namely,
for any $x,y\in A_{\theta/2,\lambda}$ we have
\begin{align*}
|\hat n_j^\rho(x)-\hat n_j^\rho(y)|^2
&
\lesssim \lambda^2 \sup_{A_{\theta/2,\lambda}}|\nabla \hat n_j^\rho|^2
%\\
%&
\lesssim  
\frac{1}{\theta^5\lambda^3 }
\int_{A_{\theta,\lambda}}
| \hat n_j^\rho - n_\infty|^2
\, dz\,.
\end{align*}
Taking $x\in \partial B_\lambda$ and
integrating with respect to $y$, this implies
\begin{align*}
\sup_{x\in \partial B_\lambda}
\bigg|\hat n_j^\rho(x) - \Xint{-}_{A_{\theta/2,\lambda}}\!\!\hat n_j^\rho(y)\, dy  
\bigg|^2
\lesssim 
\frac{1}{\theta^5\lambda^3  }
\int_{A_{\theta,\lambda}}
| \hat n_j^\rho - n_\infty|^2
\, dz
\,.
\end{align*}
Inserting $n_\infty$ gives
\begin{align*}
\sup_{\partial B_\lambda}|\hat n_j ^\rho -n_\infty|^2
&
\lesssim  \bigg| n_\infty - \Xint{-}_{A_{\theta/2,\lambda}}\!\!\hat n_j^\rho (y)\, dy\bigg|^2
 +
\frac{1}{\theta^5\lambda^3 }
\int_{A_{\theta,\lambda}}
| \hat n_j^\rho - n_\infty|^2
\, dy\,,
\\
&
\lesssim \Xint{-}_{A_{\theta/2,\lambda}}|\hat n_j^\rho -n_\infty|^2\, dy +
\frac{1}{\theta^5\lambda^3 }
\int_{A_{\theta,\lambda}}
| \hat n_j^\rho - n_\infty|^2
\, dy\,,
\end{align*}
which, for $\theta=1/4$, implies \eqref{eq:ptwise_estim_ninfty}.
\end{proof}

\subsection{Compactness of rescaled minimizers}\label{ss:comp_rescaled}

In this section we exploit another classical property: 
minimizing harmonic maps are compact in $H^1_{\loc}$  \cite{HKL88,luckhaus88}.
We apply the proof of this property to obtain
that the limit of $\hat n_j^\rho$ along a subsequence $\rho\to 0$ is a global minimizer of $\widehat E_j$.
From there, little extra effort is required to deduce the strong
$L^2$
 convergence of $\nabla n_j^\rho$ in $B_{1/\rho}\setminus\hat\omega_j$, see \eqref{eq:basic_comp},
so we include a proof of that fact, 
even though  Proposition~\ref{p:low} is only going to require the $H^1_{\loc}$ compactness.

\begin{lem}\label{l:basic_comp}
For any sequence $\rho\to 0$ 
and $j\in\{1,\ldots,N\}$, there exists a minimizer $\hat m_j$ of $\widehat E_j$ and
a subsequence, still denoted $\rho\to 0$, such that
\begin{align}\label{eq:basic_comp}
\int_{B_{1/\rho}\setminus\hat \omega_j}
|\nabla \hat n_j^\rho -\nabla\hat m_j|^2 
\, dx \longrightarrow 0,
\end{align}
as $\rho\to 0$.
\end{lem}

\begin{proof}[Proof of Lemma~\ref{l:basic_comp}]
Since $|x_i-x_j|\geq 2$ for $i\neq j$ we have
\begin{align*}
E_\rho(n_\rho)
\ \geq \
\sum_{j=1}^N E_\rho(n_\rho;B_{1}\setminus \omega_{j,\rho})
\ = \
\sum_{j=1}^N \hat E_j(\hat n_j^\rho;B_{1/\rho}\setminus \hat\omega_j)
\, .
\end{align*}
Combining this with the upper bound of Proposition~\ref{p:up}
we deduce 
\begin{align}\label{eq:sumEjmuj}
\sum_{j=1}^N 
\hat E_j(\hat n_j^\rho;B_{1/\rho}\setminus \hat\omega_j)
\leq \sum_{j=1}^N \mu_j +\mathcal O(\rho).
\end{align}
We use this bound in order to extract more information on the convergence of $\hat n_j^\rho$.

First recall that
$n_\rho$ has been extended to $\R^3$ 
so as to minimize the Dirichlet energy  inside each
hole $\omega_{\ell,\rho}$,
and so the rescaled 
map
 $\hat n_\ell^\rho$ minimizes the Dirichlet energy inside $\hat\omega_\ell$.
Moreover, we can construct an energy competitor
$w\in H^1(\hat\omega_\ell;\mathbb S^2)$
such that $w=\hat n_\ell^\rho$ on $\partial\hat\omega_\ell$
and 
\begin{align*}
\int_{\hat\omega_\ell}|\nabla w|^2\, dx \lesssim 
\int_{B_2\setminus\hat \omega_\ell}|\nabla\hat n_\ell^\rho|^2\, dx\,.
\end{align*}
This follows by applying \cite[Lemma~A.1]{HKL88} (the proof of which is valid in any domain)
to an $\R^3$-valued extension 
 with the same estimate.
The existence of this $\R^3$-valued extension follows e.g. from composing a bounded extension operator
$H^{1/2}(\partial\hat\omega_\ell)\to H^1(\hat\omega_\ell)$ 
 with the trace operator
 $H^1(B_2\setminus \hat \omega_\ell)\to H^{1/2}(\partial\hat\omega_\ell)$.
Thanks to that energy competitor $w$,
 the minimality of $\hat n_\ell^\rho$ in $\hat\omega_\ell$  implies
\begin{align*}
\int_{\hat\omega_\ell}|\nabla \hat n_\ell^\rho|^2\, dx \lesssim \int_{B_2\setminus\hat \omega_\ell}|\nabla\hat n_\ell^\rho|^2\, dx\,,\qquad\text{for }\ell=1\,\ldots,N\,.
\end{align*}
As a result, the Dirichlet energy of $\hat n_j^\rho$ in the whole space $\R^3$
is controlled by
\begin{align*}
\int_{\R^3}|\nabla\hat n_j^\rho|^2\, dx
\leq E_\rho(n_\rho) + \sum_{\ell=1}^N  \int_{\hat \omega_\ell}|\nabla\hat n_\ell^\rho|^2 dx
\lesssim E_\rho(n_\rho)\,.
\end{align*}
Using also Hardy's inequality,
we deduce
\begin{align}\label{eq:hardyhatnrhoj}
\int_{\R^3}\frac{|\hat n_j^\rho-n_\infty|^2}{1+|x|^2}\, dx 
&
\lesssim \int_{\R^3}|\nabla\hat n_j^\rho|^2\, dx
 \lesssim E_\rho(n_\rho)\,.
\end{align}
Thanks to this bound,
 for every~$j \in \{1,\ldots,N\}$,  there exists a map $\hat m_j\in H^1_{\loc}(\R^3;\mathbb S^2)$
and a subsequence $\rho\to 0$ such that
\begin{align*}
&
\nabla\hat n_j^{\rho}\rightharpoonup \nabla\hat m_j
\text{ weakly in }L^2(\R^3),
\\
&
\hat n_j^{\rho}\rightharpoonup \hat m_j
\text{ weakly in }L^2\left(\R^3;\frac{dx}{1+|x|^2}\right),
\\
\text{ and }
&\int_{\R^3} \frac{|\hat m_j-n_\infty|^2}{1+|x|^2}\, dx
\lesssim \int_{\R^3}|\nabla\hat m_j|^2\, dx <\infty.
\end{align*}
Since $\hat n_j^\rho$ minimizes the Dirichlet energy locally in $\R^3\setminus \overline{\hat \omega_j}$, 
the compactness results of \cite{HKL88,luckhaus88}
imply that $\hat n_j^\rho\to \hat m_j$ strongly in $H^1_{\loc}(\R^3\setminus \overline{\hat\omega_j})$,
and $\hat m_j$ is a local minimizer of the Dirichlet energy.
Next we adapt these arguments to show that 
 $\hat m_j$ is a global minimizer of the energy $\widehat E_j$.

To that end we fix a competitor $m\in H^1_{\loc}(\R^3\setminus \hat \omega_j;\mathbb S^2)$ such that
\begin{align*}
\int_{\R^3\setminus \hat \omega_j}\frac{ |m-n_\infty|^2}{1+|x|^2}\, dx <\infty\,,
\end{align*}
  and a radius $R>0$.
The argument in \cite[Proposition~5.1]{HKL88} provides  a sequence $\delta_\rho\searrow 0$ and 
a map $u_\rho\in H^1_{\loc}( \R^3\setminus \hat\omega_j; S^2)$
such that 
\begin{align}\label{eq:extensionHKL}
&
u_\rho=
\begin{cases}
m & \text{ in }B_R,
\\
\hat n_j^\rho
&
\text{ in }\R^3\setminus B_{R+\delta_\rho},
\end{cases}
\qquad\text{and}
\quad
\int_{B_{R+\delta_\rho}\setminus B_R}|\nabla u_\rho|^2\, dx \to 0\,.
\end{align}
The minimality of $n_\rho$ for $E_\rho$
implies $E_\rho(n_\rho)\leq E_\rho(u_\rho((\cdot-x_j)/\rho))$.
Denoting by $\Omega_j^\rho$ the rescaled domain
\begin{align*}
\Omega_j^\rho =\frac 1\rho (\Omega_\rho-x_j)\,,
\end{align*}
and taking into account the properties \eqref{eq:extensionHKL} of $u_\rho$,
this turns into
\begin{align*}
E_\rho(n_\rho)
&
=\int_{\Omega_j^\rho}|\nabla \hat n_j^\rho|^2\, dx
+F_j(\hat n_j^\rho\lfloor\partial\hat\omega_j) 
+\sum_{i\neq j} F_i(\hat n_j^\rho \lfloor\partial\hat\omega_i)
\\
&
\leq 
\int_{\Omega_j^\rho} |\nabla u_\rho|^2\, dx
+ F_j(u_\rho\lfloor \partial\hat\omega_j)
+\sum_{i\neq j} F_i ( u_\rho((x_i-x_j)/\rho +\,\cdot\,)\lfloor\partial\hat\omega_i)
\\
&
=\int_{B_R\setminus \hat\omega_j}|\nabla m|^2\, dx  
+F_j(m\lfloor\partial\hat\omega_j)
+
\int_{\Omega_j^\rho\setminus B_{R+\delta_\rho}}|\nabla \hat n_j^\rho|^2\, dx
\\
&\quad
+\sum_{i\neq j} F_i(\hat n_i^\rho \lfloor\partial\hat\omega_i)+o(1)\,,
\end{align*}
and therefore
\begin{align}\label{eq:nrhomin}
&
\int_{B_R\setminus \hat\omega_j}|\nabla \hat n_j^\rho|^2\, dx
+ F_j(\hat n_j^\rho\lfloor\partial\hat\omega_j)
\nonumber
\\
&
\leq 
\int_{B_{R+\delta_\rho}\setminus \hat\omega_j}|\nabla \hat n_j^\rho|^2\, dx
+ F_j(\hat n_j^\rho\lfloor\partial\hat\omega_j)
\nonumber
\\
&
=
E_\rho(n_\rho) - \int_{\Omega_j^\rho\setminus B_{R+\delta_\rho}}|\nabla \hat n_j^\rho|^2\, dx
-
\sum_{i\neq j} F_i(\hat n_i^\rho \lfloor\partial\hat\omega_i)
\nonumber
\\
&
\leq
\int_{B_R\setminus \hat\omega_j}|\nabla m|^2\, dx  
+F_j(m\lfloor\partial\hat\omega_j) +o(1).
\end{align}
Both terms in the 
first line
 of \eqref{eq:nrhomin}
are lower semicontinuous with respect to the weak convergence 
$\hat n_j^\rho \rightharpoonup \hat m_j$ in $H^1(B_R)$ and weak convergence of the traces in $H^{1/2}(\partial\hat\omega_j)$. 
So we deduce
\begin{align*}
\int_{B_R\setminus\hat\omega_j}|\nabla \hat m_j|^2\, dx  +F_j(\hat m_j\lfloor \partial\hat\omega_j)
\leq 
\int_{B_R\setminus \hat\omega_j}|\nabla m|^2\, dx  
+F_j(m\lfloor\partial\hat\omega_j)\,,
\end{align*}
and sending $R\to +\infty$ 
we conclude that $\hat m_j$ is a minimizer of $\widehat E_j$.

Moreover, applying the inequality \eqref{eq:nrhomin} to $m=\hat m_j$, 
and using again the lower semicontinuity of both terms in its left-hand side, 
we deduce the chain of inequalities
\begin{align*}
&
\int_{B_R\setminus \hat\omega_j}|\nabla \hat m_j|^2\, dx  
+F_j(\hat m_j\lfloor\partial\hat\omega_j)
\\
&
\leq
\liminf_{\rho\to 0} 
\int_{B_R\setminus \hat\omega_j}|\nabla \hat n_j^\rho|^2\, dx
+\liminf_{\rho\to 0}  
F_j(\hat n_j^\rho \lfloor\partial\hat\omega_j)
\\
&
\leq
\liminf_{\rho\to 0} 
\Big(
\int_{B_R\setminus \hat\omega_j}|\nabla \hat n_j^\rho|^2\, dx
+F_j(\hat n_j^\rho \lfloor\partial\hat\omega_j)
\Big)
\\
&
\leq
\limsup_{\rho\to 0} 
\Big(
\int_{B_R\setminus \hat\omega_j}|\nabla \hat n_j^\rho|^2\, dx
+F_j(\hat n_j^\rho \lfloor\partial\hat\omega_j)
\Big)
\\
&
\leq 
\int_{B_R\setminus \hat\omega_j}|\nabla \hat m_j|^2\, dx  
+F_j(\hat m_j\lfloor\partial\hat\omega_j)\,.
\end{align*}
All these inequalities must therefore be equalities, 
which implies
%\begin{align*}
%\int_{B_R\setminus \hat\omega_j}|\nabla \hat m_j|^2\, dx  
%+F_j(\hat m_j\lfloor\partial\hat\omega_j)
%&
%=\lim_{\rho\to 0}
%\Big(
%\int_{B_R\setminus \hat\omega_j}|\nabla \hat n_j^\rho|^2\, dx
%+F_j(\hat n_j^\rho \lfloor\partial\hat\omega_j) \Big)\,,
%\\
%\int_{B_R\setminus \hat\omega_j}|\nabla \hat m_j|^2\, dx 
%&
%=
%\liminf_{\rho\to 0} 
%\int_{B_R\setminus \hat\omega_j}|\nabla \hat n_j^\rho|^2\, dx
%\,,
%\\
%F_j(\hat m_j\lfloor\partial\hat\omega_j)
%&
%=
%\liminf_{\rho\to 0}  
%F_j(\hat n_j^\rho \lfloor\partial\hat\omega_j)\,,
%\end{align*}
%and eventually
\begin{align}
\label{eq:comp_BR}
%\begin{aligned}
\int_{B_R\setminus \hat\omega_j}|\nabla \hat m_j|^2\, dx 
&
=
\lim_{\rho\to 0} 
\int_{B_R\setminus \hat\omega_j}|\nabla \hat n_j^\rho|^2\, dx
\,,
\nonumber
\\
F_j(\hat m_j\lfloor\partial\hat\omega_j)
&
=
\lim_{\rho\to 0}  
F_j(\hat n_j^\rho \lfloor\partial\hat\omega_j)\,.
%\end{aligned}
\end{align}
By definition of $\mu_j=\widehat E_j(\hat m_j)$,
for any $\e>0$
we can choose $R>1$ such that
\begin{align*}
\mu_j-\e
&
\leq 
\widehat E_j(\hat m_j;B_R\setminus\hat\omega_j) 
= \lim_{\rho\to 0}
\widehat E_j(\hat n_j^\rho; B_R\setminus \hat \omega_j)\,.
\end{align*}
The last equality follows from
 \eqref{eq:comp_BR}.
 Since this is valid for any $\e>0$,
  we infer
\begin{align*}
\liminf_{\rho\to 0} \hat E_j(\hat n_j^{\rho};B_{1/\rho}\setminus\hat\omega_j) \geq \mu_j,
\end{align*}
for all $j\in\lbrace 1,\ldots,N\rbrace$.
Combining this  with \eqref{eq:sumEjmuj} implies
\begin{align*}
\hat E_j(\hat n_j^{\rho};B_{1/\rho}\setminus\hat\omega_j) \to \mu_j,\quad\text{as }\rho\to 0
\, ,
\end{align*}
and, since by  \eqref{eq:comp_BR} we also have 
$F_j(\hat n_j^\rho \lfloor\partial\hat\omega_j)
\to
F_j(\hat m_j\lfloor\partial\hat\omega_j)
$, 
\begin{align*}
\int_{B_{1/\rho}\setminus\hat\omega_j}|\nabla \hat n_j^{\rho}|^2\, dx -
\int_{B_{1/\rho}\setminus\hat\omega_j}|\nabla \hat m_j|^2\, dx \to 0.
\end{align*}
We deduce
\begin{align*}
\int_{B_{1/\rho}\!\setminus\hat\omega_j}|\nabla\hat n_j^{\rho}-\nabla\hat m_j|^2\, dx
&
=\int_{B_{1/\rho}\!\setminus\hat\omega_j}|\nabla \hat m_j|^2\, dx 
- \int_{B_{1/\rho}\!\setminus\hat\omega_j}|\nabla \hat n_j^{\rho}|^2\, dx 
\nonumber
\\
&
\quad
+ 2 \int_{\R^3\!\setminus\hat\omega_j}\langle
 \mathbf 1_{B_{1/\rho}} \nabla \hat m_j,\nabla \hat n_j^{\rho}-\nabla\hat m_j\rangle \, dx 
\nonumber
\\
&
\to 0,
\end{align*}
thanks to the weak convergence $\nabla\hat n_j^{\rho}-\nabla\hat m_j \rightharpoonup 0$ and the strong convergence $\mathbf 1_{B_{1/\rho}}\nabla \hat m_j\to\nabla\hat m_j$ in $L^2(\R^3\setminus\hat\omega_j)$.
\end{proof}

\subsection{Lower bound in terms of $\hat n_\rho^j-\hat m_j$}
\label{ss:lownjmj}

Recall that our goal is to obtain the asymptotic expansion
\begin{align*}
E_\rho(n_\rho)
&
=\sum_j \mu_i
-4\pi\rho\sum_{i\neq j}
 \frac{\langle v_i,v_j\rangle}{|x_i-x_j|}
+o(\rho)
%\\
\qquad
\text{ as }\rho\to 0.
\end{align*}
The upper bound was obtained in Proposition~\ref{p:up},
via a competitor equal to
the rescaled single-particle minimizers 
$\hat m_j$ inside small balls $B_\sigma(x_j)$,
 harmonically extended  (and projected onto $\mathbb S^2$) outside these balls.
In particular, at the gluing scale $\sigma$ around each $x_j$, that competitor was equal to the rescaled $\hat m_j$.
In this section we establish a converse estimate: 
if there is a scale $\sigma=\lambda\rho$ 
such that $n_\rho$ is close enough  to the rescaled $\hat m_j$ near $\partial B_\sigma(x_j)$, 
then the lower bound is satisfied with a small error.
Such an estimate is natural: 
 here the important point is that we manage to obtain one that is sharp enough to conclude using only the convergence \eqref{eq:basic_comp}.

\begin{prop}\label{p:low_Theta}
There exist $C>0$ and $\lambda_0\geq 2$ such that
\begin{align}
E_\rho(n_\rho)
&
\geq 
\sum_j \mu_i
-4\pi\rho\sum_{i\neq j}
 \frac{\langle v_i,v_j\rangle}{|x_i-x_j|}
% \nonumber
% \\
% &
% \quad
 -C \rho\, 
 \Xi_\lambda(\rho) \,,
 \nonumber
 \\
\text{where}
\quad
\Xi_\lambda(\rho)
&
:=
\lambda \rho 
+\frac 1\lambda
+\Theta_\lambda(\rho)^{1/2} 
+\Theta_\lambda(\rho) 
+\frac{1 +\Theta_\lambda(\rho)^2}{\lambda^3\rho}\,,
\label{eq:Theta}
\\
\text{and}
\quad
\Theta_\lambda(\rho)
&
:=\lambda^2 
\sum_j \Xint{-}_{ B_{2\lambda}\setminus B_{\lambda/2}} 
|\hat n_j^\rho -\hat m_j|^2\, dx\,,
\nonumber
\end{align}
for all
$\rho\in (0,1/(2\lambda_0))$
and $\lambda\in [\lambda_0,1/(2\rho)]$.
\end{prop}

\begin{rem}\label{r:low_Theta}
Recall that 
 $|\hat m_j-n_\infty|^2\lesssim 1/|x|^2$ for $|x|\gg 1$,
due to the asymptotic expansion \eqref{eq:vj}.
It seems  reasonable 
to hope that $|\hat n_j^\rho-n_\infty|^2$ 
could satisfy the same bound, which would imply $\Theta_\lambda(\rho)\lesssim 1$.
If we manage to
take this to the next order 
 and find a scale
$\lambda=\lambda_\rho$ such that
\begin{align*}
\Theta_{\lambda_\rho}(\rho)\ll 1\,,
\quad \lambda_\rho \ll \frac 1\rho\,,
\quad
\text{and }
\lambda_\rho \gg \frac{1}{\rho^{1/3}}\,,
\end{align*}
then we deduce that the
 error in \eqref{eq:Theta} 
satisfies $\Xi_{\lambda_\rho}(\rho)\ll 1$.
This is precisely how, in the next section, we are going to prove the lower bound of Proposition~\ref{p:low}.
\end{rem}

The proof of Proposition~\ref{p:low_Theta}
relies on two separate lower bounds:
 in the domain $U_{\lambda\rho}$ outside the balls $B_{\lambda\rho}(x_j)$, and in each
ball $B_{\lambda\rho}(x_j)$.
Specifically, we establish:
\begin{itemize}
\item 
In Lemma~\ref{l:low_ext_theta}, a lower bound for $E_\rho(n_\rho;U_{\lambda\rho})$
 in terms of the  boundary values of $\hat n_j^\rho$ at $\partial B_\lambda$.
It is simply obtained as
 the energy of the harmonic extension of $n_\rho$ from $\partial U_{\lambda\rho}$, for which
Proposition~\ref{p:harmextspheres} provides a precise expression in terms of the boundary values.
\item
In Lemma~\ref{l:low_int_theta},
a lower bound
 for 
$\widehat E_j(\hat n_j^\rho;B_\lambda)$ in terms of 
$\mu_j=\widehat E_j(\hat m_j)$
  and the boundary values of
 $\hat n_j^\rho$ at $\partial B_\lambda$.
 It follows from an upper bound on $\mu_j$ 
 obtained by constructing a competitor
equal to $\hat n_j^\rho$ inside $B_\lambda$,
and equal to its $\mathbb S^2$-projected harmonic extension outside $B_\lambda$.
 \end{itemize}
When summing these two lower bounds, 
it turns out that the main contributions from
the boundary values of
 $\hat n_j^\rho$ at $\partial B_\lambda$ 
 cancel each other,
leaving us with the rather precise lower bound
 of Proposition~\ref{p:low_Theta}.

Both  lower bounds are expressed in terms of the
spherical harmonics coefficients
\begin{align}\label{eq:hatajk}
\hat a_k^j(\lambda,\rho)
&
=\int_{\mathbb S^2}(\hat n_j^\rho(\lambda \omega)-n_\infty) \Phi_k(\omega)\, d\mathcal H^2(\omega)
\, .
\end{align}
We start with the lower bound in the exterior domain $U_{\lambda\rho}$. 

\begin{lem}\label{l:low_ext_theta}
There exist $C>0$ and $\lambda_0\geq 2$ such that
\begin{align}\label{eq:low_ext_theta}
E_\rho(n_\rho; U_{\lambda\rho})
&
\geq 
\lambda\sum_j \sum_{k\geq 0}\gamma_k^- |\hat a_k^j(\lambda,\rho)|^2
-
4\pi\rho\sum_{i\neq j}\frac{\langle v_j,v_j\rangle}{|x_i-x_j|}
\nonumber
\\
&
\quad
-C\,\rho\, \Big(
\theta_\lambda(\rho) 
+ 
\theta_\lambda(\rho)^{1/2}
+
\frac 1\lambda  
+\lambda\rho
\Big)\,,
\end{align}
where
\begin{equation}
\theta_\lambda(\rho)
:=\lambda^2\sum_j \Xint{-}_{\partial B_\lambda}|\hat n_j^\rho -\hat m_j|^2\, d\mathcal H^2\,,
\label{eq:thetalambda}
\end{equation}
for all $\rho\in (0,1/\lambda_0)$ and 
$\lambda\in [\lambda_0,1/\rho]$.
\end{lem}
Then,
complementary to the
exterior lower bound of Lemma~\ref{l:low_ext_theta},
 we have the following interior lower bound for each ball $B_{\lambda\rho}(x_j)$.

\begin{lem}\label{l:low_int_theta}
There exist $C>0$ and $\lambda_0\geq 2$ such that 
\begin{align}\label{eq:low_int_theta}
\widehat E_j(\hat n_j^\rho;B_\lambda\setminus\hat \omega_j)
%\nonumber
%\\
&
\geq
\mu_j 
- 
\lambda\sum_{k\geq 0}\gamma_k^- |\hat a_k^j(\lambda,\rho)|^2
- C\frac{1+\widetilde\Theta_\lambda(\rho)^2}{\lambda^3}\,,
\\
\text{where}
\quad
\widetilde\Theta_\lambda(\rho)
&
:=\lambda^2 \Xint{-}_{B_{5\lambda/4}\setminus B_{3\lambda/4}}
\!\!
|\hat n_j^\rho - \hat m_j|^2\, dx\,,
\nonumber
\end{align}
for all $\rho\in (0,1/(2\lambda_0))$
%$j\in\lbrace 1,\ldots,N\rbrace$ 
and 
$\lambda\in [\lambda_0,1/(2\rho)]$.
\end{lem}

Before proving the lower bounds of Lemma~\ref{l:low_ext_theta} and Lemma~\ref{l:low_int_theta},
we give the quick proof of how they imply 
Proposition~\ref{p:low_Theta}.

\begin{proof}[Proof of Proposition~\ref{p:low_Theta}]
For any $\rho\in (0,1/(2\lambda_0))$ and  $\lambda\in [\lambda_0,1/(2\rho)]$,
we can find $\lambda'\in [3\lambda/4,5\lambda/4]$ such that
$\theta_{\lambda'}(\rho)$ 
is less than its average over that interval, and then we have
\begin{align*}
\theta_{\lambda'}(\rho)\lesssim \Theta_{\lambda}(\rho),
\qquad
\widetilde \Theta_{\lambda'}(\rho)\lesssim \Theta_\lambda(\rho)\,.
\end{align*}
Summing the lower bounds \eqref{eq:low_ext_theta} and \eqref{eq:low_int_theta} taken at $\lambda=\lambda'$, we deduce the lower bound 
of Proposition~\ref{p:low_Theta}.
%\eqref{eq:Theta}.
\end{proof}

Now we prove the lower bound \eqref{eq:low_ext_theta} in the exterior domain $U_{\lambda\rho}$.

\begin{proof}[Proof of Lemma~\ref{l:low_ext_theta}]
The map $n_\rho-n_\infty$ 
 has  higher Dirichlet energy in the domain $U_{\lambda\rho}$
 than  the harmonic extension of its boundary values. 
 The harmonic extension is given by
 \begin{align*}
 n_\rho(x_j +\lambda\rho\omega)-n_\infty=
 \hat n_j^\rho(\lambda \omega)-n_\infty
 =
 \sum_{k\geq 0} 
\hat a_k^j(\lambda,\rho) \Phi_k(\omega)
 \,,
 \end{align*} 
 Proposition~\ref{p:harmextspheres} 
 provides the  lower bound
\begin{align*}
E_\rho(n_\rho;U_{\lambda\rho})
&
=\frac 1\rho \int_{U_{\lambda\rho}}|\nabla n_\rho|^2\, dx
\\
&
\geq
\lambda\sum_j \sum_{k\geq 0} \gamma_k^-|\hat a_k^j(\lambda,\rho)|^2\\
&\quad
-\lambda^2\rho\sum_j\sum_{i\neq j}\frac{
\langle \hat a_0^i(\lambda,\rho),
\hat a_0^j(\lambda,\rho)\rangle}{|x_i-x_j|}
\\
&
\quad
+
\mathcal O\left(\lambda^3\rho^2\right)\|\hat a(\lambda,\rho)\|^2_{\ell^2}\,.
\end{align*}
We define
\begin{align}\label{eq:hatb0j}
\hat b_0^j(\lambda,\rho)
&
=\hat a_0^j(\lambda,\rho)-\frac {2\sqrt\pi}{\lambda} v_j\,.
%\\
%&
%=\frac{1}{2\sqrt\pi}\int_{\mathbb S^2}\Big(\hat n_\rho^j(\lambda\omega)-n_\infty-\frac 1\lambda v_j\Big) d\mathcal H^2(\omega)\,.
%\nonumber
\end{align}
With this notation, 
we can rewrite the scalar product
$\langle \hat a_0^i(\lambda,\rho),
\hat a_0^j(\lambda,\rho)\rangle$
as
\begin{align*}
\langle \hat a_0^i(\lambda,\rho),
\hat a_0^j(\lambda,\rho)\rangle
&
=
\Big\langle
\frac {2\sqrt\pi}{\lambda} v_i
+
\hat b_0^i(\lambda,\rho),
\frac {2\sqrt\pi}{\lambda} v_j
+
\hat b_0^j(\lambda,\rho)
\Big\rangle
\\
&
=
\frac{4\pi}{\lambda^2}\langle v_i,v_j\rangle
+\mathcal O\Big( |\hat b_0(\lambda,\rho)|^2 + \frac{|\hat b_0(\lambda,\rho)|}{\lambda}\Big)\,.
\end{align*}
Hence the above lower bound becomes
\begin{align}
E_\rho(n_\rho;U_{\lambda\rho})
&
\geq 
\lambda\sum_j \sum_{k\geq 0}\gamma_k^- |\hat a_k^j(\lambda,\rho)|^2
-
4\pi\rho\sum_{i\neq j}\frac{\langle v_j,v_j\rangle}{|x_i-x_j|}
\nonumber
\\
&\quad
- C \,\rho\,
\Big(
\lambda^2|\hat b_0(\lambda,\rho)|^2 + 
\lambda | \hat b_0(\lambda,\rho)|
+\lambda^3\rho \|\hat a(\lambda,\rho)\|_{\ell^2}^2
 \Big)\,.
 \label{eq:low_ext_ab}
\end{align}
Next we  estimate the error terms in the last line.
Recalling the definitions \eqref{eq:hatb0j} of $\hat b_0^j$ and  \eqref{eq:hatajk} of $\hat a_0^j$,
 and the fact that $\Phi_0=1/(2\sqrt{\pi})$,
 we have
\begin{align*}
\hat b_0^j(\lambda,\rho)
&
=\frac{1}{2\sqrt\pi}\int_{\mathbb S^2}\Big(\hat n_\rho^j(\lambda\omega)-n_\infty-\frac 1\lambda v_j\Big)
\, d\mathcal H^2(\omega)\,.
\end{align*}
Recalling also the asymptotic expansion \eqref{eq:mj-wj} of $\hat m_j$, we can further rewrite this as
\begin{align*}
\hat b_0^j(\lambda,\rho)
&
=\frac{1}{2\sqrt\pi}\int_{\mathbb S^2}\Big(
\hat n_\rho^j(\lambda\omega)-\hat m_j(\lambda\omega)
\Big)\,d\mathcal H^2(\omega)
+\mathcal O\Big(\frac 1{\lambda^2}\Big)\,.
\end{align*}
This implies
\begin{align*}
|\hat b_0^j(\lambda,\rho)|^2
&
\lesssim \int_{\mathbb S^2}|\hat n_j^\rho-\hat m_j|^2(\lambda\omega)\, d\mathcal H^2(\omega) +\mathcal O\Big(\frac{1}{\lambda^4}\Big)
\nonumber
\\
&
\lesssim
\Xint{-}_{\partial B_\lambda}|\hat n_j^\rho-\hat m_j|^2\, d\mathcal H^2 + \mathcal O\Big(\frac{1}{\lambda^4}\Big)\,\,.
\end{align*}
Hence, recalling the definition \eqref{eq:thetalambda} of $\theta_\lambda$,
\begin{align}\label{eq:estim_hatb0j}
\lambda^2|\hat b_0(\lambda,\rho)|^2
\lesssim\theta_\lambda(\rho) +\frac{1}{\lambda^2}\,.
\end{align}
Moreover, by definition \eqref{eq:hatajk} of the coefficients $\hat a_k^j$ and by orthonormality of $(\Phi_k)$ in $L^2(\mathbb S^2)$ we have
\begin{align*}
\|\hat{a}^j(\lambda,\rho)\|_{\ell^2}^2
&
=\sum_{k\geq 0}|\hat{a}_k^j(\lambda,\rho)|^2
=\int_{\mathbb S^2}|\hat n_j^\rho(\lambda\omega)-n_\infty|^2\,d\mathcal H^2(\omega)
\\
&
\lesssim \Xint{-}_{\partial B_\lambda}
|\hat n_j^\rho -n_\infty|^2\,d\mathcal H^2
\\
&
\lesssim 
\Xint{-}_{\partial B_\lambda}
|\hat n_j^\rho -\hat m_j|^2\,d\mathcal H^2
+\sup_{\partial B_\lambda}|\hat m_j -n_\infty|^2\,.
\end{align*}
Recalling the asymptotic expansion \eqref{eq:mj-wj} of $\hat m_j$, we are left with
\begin{align*}
\lambda^2\|a^j(\lambda,\rho)\|_{\ell^2}^2
&
\lesssim
1+
\lambda^2\Xint{-}_{\partial B_\lambda}
|\hat n_j^\rho -\hat m_j|^2\,d\mathcal H^2
%\\
%&
\lesssim 1+\theta_\lambda(\rho)\,,
\end{align*}
where the last inequality follows from the definition \eqref{eq:thetalambda} of $\theta_\lambda$.
Combining this with the bound \eqref{eq:estim_hatb0j} on $\hat b_0$  we deduce
\begin{align*}
\lambda^2|\hat b_0(\lambda,\rho)|^2 + 
\lambda |\hat b_0(\lambda,\rho)|
+\lambda^3\rho \|\hat a(\lambda,\rho)\|_{\ell^2}^2
\lesssim
\theta_\lambda(\rho) 
+ 
\theta_\lambda(\rho)^{1/2}
+
\frac 1\lambda  
+\lambda\rho\,.
\end{align*}
Plugging this into \eqref{eq:low_ext_ab} we obtain \eqref{eq:low_ext_theta}.
\end{proof}

%%%%%%%%%%%%%%%%%%%%%%%%%%%%%%%%%%%%%%%%%%%%%%%%%%
%%%%%%%%%%%%%%%%%%%%%%%%%%%%%%%%%%%%%%%%%%%%%%%%%%
%%%%%%%%%%%%%%%%%%%%%%%%%%%%%%%%%%%%%%%%%%%%%%%%%%

Finally we prove the lower bound \eqref{eq:low_int_theta} in each ball $B_{\lambda\rho}(x_j)$.

\medskip

\begin{proof}[Proof of Lemma~\ref{l:low_int_theta}]
First note that, if $\lambda_0$ is large enough, 
then thanks to the pointwise estimate \eqref{eq:ptwise_estim_ninfty}
and Hardy's inequality \eqref{eq:hardyhatnrhoj}, we have 
\begin{align*}
\sup_{\partial B_\lambda}
|\hat n_j^\rho -n_\infty|^2
\leq \frac 12\,,
\qquad\forall \lambda\in [\lambda_0,1/\rho]\,.
\end{align*}
Now, in order to bound the minimal energy $\mu_j=\widehat E_j(\hat m_j)$ from above, we
consider a competitor $\tilde n_\rho^j\colon \R^3\to\mathbb S^2$ defined by
\begin{align*}
\tilde n_\rho^j
=
\begin{cases}
\hat n_j^\rho
&\quad\text{in }B_\lambda,
\\
\frac{n_\infty +\tilde u}{|n_\infty +\tilde u|}
&\quad\text{outside }B_\lambda,
\end{cases}
\end{align*}
where $\tilde u\colon\R^3\setminus B_\lambda\to\R^3$ is the  harmonic extension agreeing with
 $\hat n_\rho^j - n_\infty$ on $\partial B_\lambda$.
 By definition \eqref{eq:hatajk} of the coefficients $\hat a_k^j$,
the extension $\tilde u$ is given by
\begin{align*}
\tilde u(r\omega) =\sum_{k\geq 0} a_k^j(\lambda,\rho) \Big(\frac r\lambda\Big)^{\gamma_k^-}\Phi_k(\omega)\,,
\end{align*}
and its energy by
\begin{align}\label{eq:en_tildeu}
\int_{\R^3\setminus B_\lambda}|\nabla\tilde u|^2\, dx
=-\int_{\partial B_\lambda} \langle\tilde u ,\partial_r\tilde u\rangle\, d\mathcal H^2 =\lambda\sum_{k\geq 0}\gamma_k^- |a_k^j(\lambda,\rho)|^2\,.
\end{align}
By minimality of $\hat m_j$  we have
\begin{align}
\widehat E_j(\hat m_j;\R^3\setminus \hat\omega_j)
&
\leq \widehat E_j(\tilde n_\rho^j;\R^3\setminus\hat \omega_j)
\nonumber
\\
&
= \widehat E_j(\hat n_j^\rho;B_\lambda\setminus \hat\omega_j)
+\int_{\R^3\setminus B_\lambda} |\nabla \tilde n_\rho^j|^2\, dx
\nonumber
\\
&
\leq
\widehat E_j(\hat n_j^\rho;B_\lambda\setminus \hat\omega_j)
+\int_{\R^3\setminus B_\lambda} \frac{|\nabla\tilde u|^2}{|n_\infty +\tilde u|^2}\, dx\,.
\label{eq:mujtilden}
\end{align}
The last inequality 
follows from the 
inequality $|\nabla (v/|v|)|^2\leq |\nabla v|^2/|v|^2$
applied to $v=n_\infty +\tilde u$,
see \eqref{eq:norm_grad_proj}.
Since the harmonic function $\langle n_\infty,\tilde u\rangle$ is either positive or attains its minimum at the boundary $\partial B_\lambda$, we have
\begin{align*}
|n_\infty +\tilde u|^2 
&
= 1 +2 \langle n_\infty,\tilde u\rangle +|\tilde u|^2
\geq 1 - 2 \sup_{\partial B_\lambda} |\langle n_\infty,
\hat n_j^\rho-n_\infty\rangle |.
\end{align*}
Using also that
\begin{align*}
|\langle n_\infty,
\hat n_j^\rho - n_\infty\rangle|
=1-\langle n_\infty,
\hat n_j^\rho\rangle
=\frac 12 |n_\infty-
\hat n_j^\rho|^2,
\end{align*}
we deduce
\begin{align*}
|n_\infty +\tilde u|^2\geq 1- 
\sup_{\partial B_\lambda} 
|\hat n_j^\rho- n_\infty |^2
\qquad\text{in }\R^3\setminus B_\lambda\,.
\end{align*}
Since we have $|\hat n_j^\rho -n_\infty|^2\leq 1/2$ on $\partial B_\lambda$,
this implies
\begin{align*}
\frac{1}{|n_\infty +\tilde u|^2}
\leq 1 + 2 \sup_{\partial B_\lambda} 
| \hat n_j^\rho- n_\infty|^2
\qquad\text{in }\R^3\setminus B_\lambda\,.
\end{align*}
Using this to bound the last term in  the energy estimate
\eqref{eq:mujtilden} we obtain
\begin{align*}
\mu_j
&=
\widehat E_j(\hat m_j;\R^3\setminus \hat\omega_j)
\\
&
\leq \widehat E_j(\hat n_j^\rho;B_\lambda\setminus \hat\omega_j)
+\Big(1
+2\sup_{\partial B_\lambda}|\hat n_j^\rho-n_\infty|^2
\Big)
\int_{\R^3\setminus B_\lambda} |\nabla \tilde u|^2\, dx.
\\
&
=
\widehat E_j(\hat n_j^\rho;B_\lambda\setminus \hat\omega_j)
+\Big(1
+
2\sup_{\partial B_\lambda}
|\hat n_j^\rho-n_\infty|^2
\Big)
\lambda\sum_{k\geq 0}\gamma_k^- |\hat a_k^j(\lambda,\rho)|^2\,.
\end{align*}
The last equality follows from the explicit expression
\eqref{eq:en_tildeu} of the energy of $\tilde u$.
Rearranging, we deduce
\begin{align}\label{eq:low_int_ajk}
\widehat E_j(\hat n_j^\rho;B_\lambda\setminus \hat\omega_j)
&
\geq
\mu_j -\lambda\sum_{k\geq 0}\gamma_k^- 
|\hat a_k^j(\lambda,\rho)|^2 
\nonumber
\\
&
\quad
- 2 \lambda \sup_{\partial B_\lambda}
|\hat n_j^\rho-n_\infty|^2
\sum_{k\geq 0}\gamma_k^- 
|\hat a_k^j(\lambda,\rho)|^2\,.
\end{align}
Using that $\gamma_k^-\leq 1+(\gamma_k^-)^2\lesssim 1+\lambda_k$
and the definition \eqref{eq:hatajk} of the coefficients $\hat a_j^k$
we obtain
\begin{align*}
\sum_{k\geq 0}
\gamma_k^-|a_k^j(\lambda,\rho)|^2
&
\lesssim \sum_{k\geq 0} |a_k^j(\lambda,\rho)|^2
+\sum_{k\geq 0}\lambda_k |a_k^j(\lambda,\rho)|^2\\
&
=\int_{\mathbb S^2}
|\hat n_j^\rho(\lambda\omega)-n_\infty|^2\,d\mathcal H^2(\omega)
+\int_{\mathbb S^2}|\nabla_\omega\hat n_j^\rho(\lambda\omega)|^2\,
d\mathcal H^2(\omega)
\\
&
\lesssim 
\Xint{-}_{B_{5\lambda/4}\setminus B_{3\lambda/4}}
\!\!
|\hat n_j^\rho -n_\infty|^2\, dx\,.
\end{align*}
The last inequality follows from the pointwise estimates of Lemma~\ref{l:ptwise_estim}, which also imply
\begin{align*}
\sup_{\partial B_\lambda}|\hat n_j^\rho -n_\infty|^2
\lesssim
\Xint{-}_{B_{5\lambda/4}\setminus B_{3\lambda/4}}
\!\!
|\hat n_j^\rho -n_\infty|^2\, dx\,.
\end{align*}
Using the last two estimates to bound the last line in the lower bound \eqref{eq:low_int_ajk}, we deduce
\begin{align*}
\widehat E_j(\hat n_j^\rho;B_\lambda\setminus \hat\omega_j)
&
\geq
\mu_j -\lambda\sum_{k\geq 0}\gamma_k^- 
|\hat a_k^j(\lambda,\rho)|^2 
\nonumber
\\
&
\quad
- \frac{C}{\lambda^3}\Big(\lambda^2 
\Xint{-}_{B_{5\lambda/4}\setminus B_{3\lambda/4}}
\!\!
|\hat n_j^\rho -n_\infty|^2\, dx
\Big)^2\,.
\end{align*}
To conclude we recall that the asymptotics \eqref{eq:mj-wj}
of $\hat m_j$  ensure
\begin{align*}
\lambda^2\Xint{-}_{B_{5\lambda/4}\setminus B_{3\lambda/4}}
\!\!|\hat n_j^\rho -n_\infty|^2\, dx
\lesssim 1 + \widetilde \Theta_\lambda(\rho),
\end{align*}
with $\widetilde\Theta_\lambda$ as in \eqref{eq:low_int_theta}.
\end{proof}

\subsection{Far field asymptotics of $\hat n_j^\rho$}
\label{ss:farfieldhatnrho}

As noted in Remark~\ref{r:low_Theta}, 
the proof of the sharp lower bound 
of Proposition~\ref{p:low} 
relies on proving that $|\hat n_j^\rho -\hat m_j|$ on some large annulus $B_{2\lambda}\setminus B_{\lambda/2}$ is much smaller than the leading asymptotics of 
$\hat m_j - n_\infty$
which is of order $1/\lambda$.
In this section we show that 
$\hat n_j^\rho$ has an asymptotic 
expansion
similar to that of  $\hat m_j$ in \eqref{eq:mj-wj}. 
This
 will allow us to control the error terms~$\Theta_\lambda(\rho)$
and $\Xi(\rho)$ 
in \eqref{eq:Theta}, leading to the proof of Proposition~\ref{p:low}.

\begin{prop}\label{p:farfieldhatnrho}
There exist $\lambda_0>2$, $\rho_0\in (0,1)$ and, for every $\rho\in (0,\rho_0)$,
vectors $n_\infty^\rho,v_j^\rho\in\R^3$ such that
\begin{align}
&
\hat n_j^\rho(x) 
= n_\infty^\rho + \frac{v_j^\rho}{| x|} + w_j^\rho( x)\,,
\qquad
|n_\infty^\rho - n_\infty|^2 \lesssim \rho\ln^2\!\rho\,,
\label{eq:farfieldhatnrho}
\\
&
| w_j^\rho|^2
\lesssim 
\rho 
+\frac{\ln^6\! |x|}{|x|^4}\qquad\text{for }\lambda_0\leq | x|\leq  \frac{|\ln\rho|}{\sqrt\rho}\,,
\nonumber
\end{align}
and
$v_j^\rho \to v_j$
along the sequence $\rho\to 0$ provided by Lemma~\ref{l:basic_comp}.
\end{prop}

\medskip

Before proving Proposition~\ref{p:farfieldhatnrho} 
we give the short argument of how to combine it with Proposition~\ref{p:low_Theta} to
deduce the sharp lower bound on $E_\rho(n_\rho)$.

\begin{proof}[Proof of Proposition~\ref{p:low}]
Using the asymptotic expansions \eqref{eq:farfieldhatnrho} of $\hat n_j^\rho$ and \eqref{eq:mj-wj} of $\hat m_j$,
we have, for $\lambda_0\leq |x|\leq |\ln\rho|/\sqrt\rho$, 
\begin{align*}
|\hat n_j^\rho -\hat m_j|^2
\lesssim \rho\ln^2\!\rho  +\frac{\ln^6\!\rho}{|x|^4} +\frac{|v_j-v_j^\rho|^2}{|x|^2}\,,
\end{align*}
and recalling the definition of $\Theta_\lambda$ in \eqref{eq:Theta} we infer
\begin{align*}
\Theta_\lambda(\rho)
\lesssim \lambda^2\rho\ln^2\!\rho +\frac{\ln^6\!\rho}{\lambda^2} +|v_j-v_j^\rho|^2
\,,
\qquad
\text{for }
2\lambda_0\leq\lambda\leq \frac{|\ln\rho|}{2\sqrt\rho}
\,.
\end{align*}
Choosing, for small $\rho$, the admissible value
\begin{align*}
\lambda=\lambda_\rho =\frac{|\ln\rho|}{\rho^{1/4}}\,,
\end{align*}
we deduce
\begin{align*}
\Theta_{\lambda_\rho}(\rho)\lesssim  \sqrt\rho \ln^4\!\rho
 +|v_j^\rho -v_j|^2 \longrightarrow 0,
\end{align*}
along the sequence $\rho\to 0$ provided by Lemma~\ref{l:basic_comp}, 
since $v_j^\rho \to v_j$ by Proposition~\ref{p:farfieldhatnrho}.
Since we also have $\rho^{-1/3}\ll\lambda_\rho\ll\rho^{-1}$
we conclude, see Remark~\ref{r:low_Theta}, that the lower bound error $\Xi_\lambda$ in \eqref{eq:Theta} satisfies $\Xi_{\lambda_\rho}(\rho)\to 0$, thus proving Proposition~\ref{p:low}.
\end{proof}

Now we turn to the proof
 of Proposition~\ref{p:farfieldhatnrho}.
 It
is based on the strategy in 
\cite{ABLV23} for the asymptotic expansion \eqref{eq:vj} of $\hat m_j$.
That strategy relies on repeated application of two basic principles:
\begin{itemize}
\item if the Laplacian $\Delta n$ is small in some region, then $n$ is close to a classical harmonic function $u$, that is, $\Delta u=0$;
\item harmonic functions $u$ in $\R^3\setminus B_\lambda$ 
have an asymptotic expansion determined by their spherical harmonics decomposition.
\end{itemize}
Here our main issue is the last point: 
we can only hope to control $\hat n_j^\rho$ 
in an annulus $B_{1/\rho}\setminus B_\lambda$,
where the spherical harmonics decomposition can have 
radially increasing modes.
We have to take into account additional error terms coming from these increasing modes,
and this is reflected here in the fact that we are only able to control the error $\hat w_j$ in a smaller annulus, of amplitude slightly larger than $1/\sqrt\rho$.
%As we will see, we also need to be quite precise (and deal with borderline cases)
As we will see, we also need to 
deal with with borderline cases 
in our application of the first principle, 
i.e.\ that a function with small Laplacian is close to a harmonic function.

\begin{proof}[Proof of Proposition~\ref{p:farfieldhatnrho}]
We follow essentially the first two steps of \cite[Theorem~1.1]{ABLV23}, 
with
 adaptations for estimates in an annulus
$B_{1/\rho}\setminus B_\lambda$
 instead of the whole exterior domain $\R^3\setminus B_\lambda$.
 
The initial decay of $|\nabla \hat n_j^\rho|$ that we start with
is provided by the small energy estimate  
 \eqref{eq:ptwise_estim_grad1} and the fact that
 the energy $\widehat E_j(\hat n_j^\rho)$ is bounded:
 we have
 \begin{align*}
 |\nabla \hat n_j^\rho|^2\lesssim \frac{1}{|x|^3}\qquad\text{for }\lambda_0\leq |x| \leq \frac 1{2\rho}\,.
 \end{align*}
 Together with the harmonic map equation
 \begin{align*}
 -\Delta \hat n_j^\rho =|\nabla \hat n_j^\rho|^2\hat n_j^\rho\,,
 \end{align*}
this implies $|\Delta \hat n_j^\rho|\lesssim 1/|x|^3$, which is  not precise enough to capture the first decaying harmonic term 
of order $1/|x|$ in the expansion of $\hat n_j^\rho$.
 
 In order to obtain a stronger estimate on $|\nabla\hat n_j^\rho|$,
we proceed as in the alternative proof of Step~1 in \cite[Theorem~1.1]{ABLV23}
 and consider the map $g=\partial_\alpha n_j^\rho$, 
 which is pointwise orthogonal to $\hat n_j^\rho$ and solves the linearized equation
 \begin{align}\label{eq:lin_g}
 -\Delta g = 2\langle \nabla \hat n_j^\rho,\nabla g\rangle +|\nabla \hat n_j^\rho|^2 g.
 \end{align}
For $R\in [2\lambda_0, 1/(6\rho)]$ we multiply this with $\chi^2 g$ for a cut off function $\chi$ satisfying
\begin{align*}
\mathbf 1_{R\leq |\hat x|\leq 2R}\leq \chi(\hat x)\leq \mathbf 1_{R/2\leq |\hat x| \leq 3R}
\quad\text{and}
\quad
|\nabla \chi|\lesssim 1/R\,.
\end{align*}
Since $g$ is orthogonal to $\hat n_j^\rho$, the first term in the right-hand side drops out and we are left with
\begin{align*}
-\chi^2 \langle g,\Delta g\rangle = \chi^2 |\nabla \hat n_j^\rho|^2 |g|^2\,.
\end{align*}
Integrating by parts in the left-hand side, we deduce
\begin{align*}
&\int_{B_{3R}\setminus B_{R/2}}\!\! |\nabla g|^2 \chi^2 \, dx
\\
&
=
\int_{B_{3R}\setminus B_{R/2}}\!\! |\nabla \hat n_j^\rho|^2 |g|^2\chi^2\, dx
-2
\int_{B_{3R}\setminus B_{R/2}}\!\! 
\langle g, (\nabla\chi\cdot\nabla)g \rangle \chi\, dx
\\
&
\leq
\int_{B_{3R}\setminus B_{R/2}}\!\! |\nabla \hat n_j^\rho|^2\chi^2\, dx
+ 2\int_{B_{3R}\setminus B_{R/2}}\!\! |g|^2 |\nabla\chi|^2\, dx
\\
&
\quad +\frac 12 \int_{B_{3R}\setminus B_{R/2}}\!\! |\nabla g|^2 \chi^2 \, dx\,.
\end{align*}
Absorbing the last term into the left-hand side and using that
$|\nabla\chi|\lesssim 1/R$ and  
 $|g|^2\leq |\nabla\hat n_j^\rho|^2\lesssim 1/ R^3$, we infer
\begin{align*}
 \int_{B_{3R}\setminus B_{R/2}}\!\!
|\nabla g|^2\chi^2\, dx \lesssim\frac 1{R^2}\,,
\end{align*}
and therefore
\begin{align*}
\Xint{-}_{B_{2R}\setminus B_R}\! |\nabla g|^2\, dx \lesssim \frac 1{R^5}\,.
\end{align*}
Using this, and once more  $|g|^2\leq |\nabla\hat n_j^\rho|^2\lesssim 1/ R^3$, to estimate the right-hand side of \eqref{eq:lin_g}, we find
\begin{align*}
\bigg(
\Xint{-}_{B_{2R}\setminus B_R}\! |\Delta g|^2\, dx 
\bigg)^{1/2} 
\lesssim \frac 1{R^4}
\qquad
\text{for } \lambda_0\leq R \leq \frac{1}{6\rho}\,.
\end{align*}
Applying Lemma~\ref{l:decaylog} with $d=3$, $\gamma=2$ and $f=\mathbf 1_{B_{1/(6\rho)}}\Delta g$, we obtain the existence of 
a map $u\colon B_{1/(6\rho)}\setminus B_{\lambda_0}\to\R^3$
such that 
$\Delta (u-g)=0$ and
\begin{align}\label{eq:decay_harm_corr_g}
\Xint{-}_{B_{2R}\setminus B_R}\! |u|^2\, dx
\lesssim \frac{\ln^2 R}{R^4}
\qquad
\text{for } \lambda_0\leq R \leq \frac{1}{6\rho}\,.
\end{align}
This implies in particular
\begin{align*}
\int_{B_{1/(6\rho)}\setminus B_{\lambda_0}}|u|^2\, dx \lesssim 1\,.
\end{align*}
This, together with
the inequality $|g|^2\leq |\nabla \hat n_j^\rho|^2$ 
and
the fact that $\widehat E_j(\hat n_j^\rho)\lesssim 1$,
implies
\begin{align*}
\int_{B_{1/(6\rho)}\setminus B_{\lambda_0}}|u-g|^2\, dx \lesssim 1\,.
\end{align*}
We may therefore apply Lemma~\ref{l:harm_L2_annulus} to the harmonic function $u-g$.
This gives
 \begin{align*}
\Xint{-}_{B_{2R}\setminus B_R}|u-g|^2\, d x 
&
\lesssim \frac{\lambda_0}{R^4}
+
\frac{\rho}{R^2}
\qquad \text{for }2\lambda_0\leq R\leq \frac{1}{24\rho}.
\end{align*}
Combining this with the decay
\eqref{eq:decay_harm_corr_g} of $u$
and raising the value of $\lambda_0$,
 we 
deduce
 \begin{align*}
 \Xint{-}_{B_{2R}\setminus B_R}|g|^2\, dx 
&
\lesssim \frac{\ln^2 R}{R^4}
+
\frac{\rho}{R^2}
\qquad \text{for }\lambda_0\leq R\leq \frac{1}{24\rho}.
\end{align*}
Recalling the definition $g=\partial_\alpha\hat n_j^\rho$, applying this for all $\alpha=1,2,3$ and using the small energy estimate \eqref{eq:ptwise_estim_grad1}, 
we obtain
\begin{align}\label{eq:firstdecaygradnrho}
|\nabla \hat n_j^\rho|^2
\lesssim
 \frac{\ln^2 |x|}{|x|^4}
+
\frac{\rho}{|x|^2}
\qquad \text{for }\lambda_0\leq |x|\leq \frac{1}{24\rho}.
\end{align}
Since 
 $r\mapsto r^{-2}\ln^2 r$ is decreasing for $r\geq e$,
we have $\rho\lesssim r^{-2} \ln^2r$ for all $r\in [e,8|\ln\rho|/\sqrt\rho]$, 
and
 we deduce
\begin{align*}
|\Delta \hat n_j^\rho| = |\nabla \hat n_j^\rho|^2
\lesssim
 \frac{\ln^2 |x|}{|x|^4}
 \qquad\text{for }\lambda_0\leq |x|\leq R_\rho :=\frac{8|\ln\rho|}{\sqrt\rho}\,.
\end{align*}
Applying Lemma~\ref{l:decaylog} with $\gamma=2=\theta$
(and elliptic estimates to turn its conclusion into a pointwise bound),
 we find
$\tilde u\colon B_{R_\rho}\setminus B_{\lambda_0} \to \R^3$
such that $\Delta(\hat n_\rho^j-\tilde u)=0$ and
\begin{align*}
\frac{|\tilde u|}{| x|} +|\nabla \tilde u| 
\lesssim \frac{\ln^3\! | x|}{|x|^3}
\qquad \text{for }\lambda_0\leq | x| \leq
R_{\rho}\,.
\end{align*}
Since $\hat n_\rho^j -\tilde u$ is a harmonic function which satisfies
\begin{align*}
|\nabla (\hat n_\rho^j-\tilde u)|
\lesssim  \frac{\ln |x|}{|x|^2}
\qquad \text{for }\lambda_0\leq | x| \leq
R_{\rho}\,,
\end{align*}
Lemma~\ref{l:harm_H1_annulus} 
 allows us to 
decompose it as
\begin{align*}
\hat n_\rho^j -\tilde u=n_\infty^\rho
+ \frac{v_j^\rho}{|x|} +  v+\tilde w,
\end{align*} 
where
$n_\infty^\rho, v_j^\rho\in\R^3$,
 $v$ is harmonic in $\R^3\setminus B_{\lambda_0}$,
and
(possibly raising the value of $\lambda_0$),
 \begin{align*}
 &|v_j^\rho|\lesssim 1,
 \qquad
|v|+|x|\,|\nabla v| 
\lesssim \frac{1}{|x|^2}\qquad\text{for }|x|\geq \lambda_0\,,
\\
&
|\tilde w| +|x|\, |\nabla\tilde w|
\lesssim \frac{\ln R_\rho}{R_\rho}
\lesssim \sqrt\rho
\qquad \text{for }\lambda_0\leq |x|\leq \frac{R_\rho}{8}
 =\frac{|\ln\rho|}{\sqrt\rho}\,.
 \end{align*}
Setting $w_j^\rho = \tilde u +v + \tilde w$, we obtain
\begin{align}
&
\hat n_j^\rho = n_\infty^\rho +\frac{v_j^\rho}{|x|} 
+w_j^\rho
\,,
\label{eq:farfieldharnrho1}
\\
&
|w_j^\rho| +|x|\, |\nabla w_j^\rho| \lesssim
\frac{\ln^3\!|x|}{|x|^2} +\sqrt\rho
\qquad
 \text{for }\lambda_0\leq |x|\leq \frac{|\ln\rho|}{\sqrt\rho}
\,.
\nonumber
\end{align}

\medskip

To complete the proof of Proposition~\ref{p:farfieldhatnrho}, 
it remains to 
obtain the estimate 
\eqref{eq:farfieldhatnrho} 
on $|n_\infty^\rho-n_\infty|$ and 
that $|v_j^\rho - v_j|\to 0$ 
as $\rho\to 0$.
 
\medskip 
 
From the expansion \eqref{eq:farfieldharnrho1}
and the fact that $|v_j^\rho|\lesssim |\ln\rho|$, we infer
\begin{align}\label{eq:ninftyrho1}
&
|n_\infty^\rho - n_\infty|^2
=\Big| \hat n_j^\rho - n_\infty -\frac{v_j^\rho}{|x|}-w_j^\rho\Big|^2
\nonumber
\\
&
\lesssim  |\hat n_j^\rho - n_\infty|^2 
+\frac{\ln^2\!\rho}{|x|^2} + \frac{\ln^6\!\rho}{|x|^4} +\rho
\qquad
 \text{for }\lambda_0\leq |x|\leq\frac{|\ln\rho|}{\sqrt\rho}
 \,.
\end{align}
Moreover, 
the pointwise bound \eqref{eq:firstdecaygradnrho} on $|\nabla \hat n_j^\rho|$
and
the fundamental theorem of calculus
ensure,
 for $|\ln\rho|/\sqrt\rho\leq |x|\leq 1/(24\rho)$,
\begin{align*}
\sup_{\partial B_{|\ln\rho|/\sqrt\rho}}
|n_j^\rho-n_\infty|^2
&
\lesssim |\hat n_j^\rho -n_\infty|^2(x) 
+\Big(\int_{|\ln\rho|/\sqrt\rho}^{|x|}
\frac{\ln r + \sqrt\rho\, r}{r^2} 
dr
\Big)^2
\\
&
\lesssim 
|\hat n_j^\rho -n_\infty|^2(x) +\rho \ln^2 |x|\,.
\end{align*}
So the inequality \eqref{eq:ninftyrho1} at $|x|=|\ln\rho|/\sqrt\rho$ implies
\begin{align*}
|n_\infty^\rho -n_\infty|^2
\lesssim |\hat n_j^\rho -n_\infty|^2 + 
\rho\ln^2 |x|\qquad\text{ for }
\frac{|\ln\rho|}{\sqrt\rho}\leq |x|\leq \frac{1}{24\rho}\,.
\end{align*}
%Combining this with \eqref{eq:ninftyrho1} we deduce
%\begin{align*}
%|n_\infty^\rho -n_\infty|^2
%\lesssim |\hat n_j^\rho -n_\infty|^2 
%+\frac{\ln^2\!\rho}{|x|^2} + \frac{\ln^6\!\rho}{|x|^4} 
%+ 
%\rho\ln^2 |x|\,,
%\end{align*}
%in the whole range $\lambda_0\leq |x|\leq 1/(48\rho)$.
Dividing by $|x|^2$ and integrating on the 
annulus $1/(48\rho)\leq |x|\leq 1/(24\rho)$, we deduce
\begin{align*}
\frac{|n_\infty^\rho -n_\infty|^2}{\rho}
&
\lesssim \int_{B_{\frac{1}{48\rho}}\setminus 
B_{\frac{1}{24\rho}}}\!\! 
\frac{|n_\infty^\rho -n_\infty|^2}{|x|^2}\, dx 
\\
&
\lesssim \int_{B_{1/\rho}\setminus B_1}
\frac{|\hat n_j^\rho -n_\infty|^2}{|x|^2}\, dx 
+\ln^2\!\rho\,.
\end{align*}
Recalling Hardy's inequality \eqref{eq:hardyhatnrhoj}, we obtain the claimed estimate 
\begin{align*}
|n_\infty^\rho -n_ \infty|^2\lesssim \rho\ln^2\!\rho\,.
\end{align*}
Finally we turn to the estimate on $|v_j^\rho -v_j|$.
Using the expansions \eqref{eq:farfieldharnrho1} and \eqref{eq:mj-wj} of $\hat n_j^\rho$ and $\hat m_j$ we express 
\begin{align*}
%&
\frac{|v_j^\rho - v_j|^2}{|x|^4}
&
=
\bigg|\nabla \Big( \frac{v_j^\rho}{|x|}-\frac{v_j}{|x|}\Big)
\bigg|^2
%\\
%&
=
\big|\nabla (\hat n_j^\rho -
 \hat m_j 
+\hat w_j-w_j^\rho
)
\big|^2
\\
&
\lesssim  |\nabla \hat n_j^\rho -\nabla \hat m_j|^2
+ \frac{\ln^6\!|x|}{|x|^6}+ \frac{\rho}{|x|^2}
\qquad\text{for }\lambda_0\leq |x|\leq\frac{|\ln\rho|}{\sqrt\rho}\,.
\end{align*}
Integrating this inequality over an annulus 
$\lambda\leq |x|\leq 2\lambda$  for 
any $\lambda\in [\lambda_0,1/\sqrt\rho]$,
we find
\begin{align*}
\frac{|v_j^\rho-v_j|^2}{\lambda}
\lesssim 
\int_{B_{2\lambda}\setminus B_{\lambda}}
\!\!
|\nabla \hat n_j^\rho -\nabla \hat m_j|^2
\, dx
+\frac{\ln^6\!\lambda}{\lambda^3} 
+
\lambda\,\rho \,.
\end{align*}
Along the sequence $\rho \to 0$ provided by Lemma~\ref{l:basic_comp},
 the first integral in the right-hand side 
 converges to zero.
 Hence we deduce, along that sequence,
 \begin{align*}
 \limsup_{\rho\to 0} |v_j^\rho-v_j|^2
 \lesssim 
 \frac{\ln^6\!\lambda}{\lambda^2} \qquad\forall \lambda\geq \lambda_0.
 \end{align*}
 Sending $\lambda\to \infty$ concludes the proof that $v_j^\rho\to v_j$.
\end{proof}

\begin{appendices}

\section{Decaying solutions of Poisson's equation}\label{a:poisson}

We include here, for the readers' convenience, 
a proof of a folklore result
about existence of decaying solutions to Poisson's equations.
We follow and adapt the proof in \cite[Lemma~A.2]{ABLV23}
in the case $\theta=0$ and $\gamma$ non-integer.

\begin{lem}\label{l:decaylog}
Let $d\geq 3$,  $\gamma\geq d-2$ and $\theta\geq 0$,
$\lambda\geq 1$   and  $f$ a function in $\R^d\setminus B_\lambda$ satisfying 
\begin{align*} 
\left(\Xint{-}_{R<|x|<2R}  f^2 \, dx  \right)^{\frac 12}&\leq \frac{\ln^\theta \! (2 R/\lambda)}{R^{\gamma+2}}\qquad\forall R\geq \lambda.
\end{align*}
Then there exists a function $u$ such that $\Delta u =f$ in $\R^d\setminus B_\lambda$ and
\begin{align*} 
\left(\Xint{-}_{R<|x|<2R}
u^2\, dx\right)^{\frac 12}
 &\leq C \frac{\ln^{1+\theta} \!(2R/\lambda)}{R^{\gamma}}\qquad\forall R\geq \lambda,
\end{align*}
where $C>0$ depends only on $d$,  $\gamma$ and $\theta$.
\end{lem}

\begin{proof}[Proof of Lemma~\ref{l:decaylog}]
By scaling, we assume without loss of generality that $\lambda=1$.
We fix,
as in \S~\ref{s:harm_ext},
an orthonormal Hilbert basis $\lbrace \Phi_j\rbrace$  of $L^2(\mathbb S^{d-1})$ which diagonalizes the Laplace-Beltrami operator,
\begin{equation*}
-\Delta_{\mathbb S^{d-1}}\Phi_j =\lambda_j \Phi_j,\qquad 0= \lambda_0 \leq \lambda_1\leq\cdots
\end{equation*}
The set $\lbrace \lambda_j\rbrace_{j\in\mathbb N}$ coincides with $\lbrace k^2 + k(d-2)\rbrace_{k\in\mathbb N}$. 
The eigenfunctions corresponding to $k^2 + k(d-2)$ span the homogeneous harmonic polynomials 
of degree $k$. 
For  a $W^{2,2}_{loc}$ function $w\colon (0,\infty) \to \R$ we have
\begin{equation} \label{e.Lj}
\Delta (w(r)\Phi_j(\omega)) =(\mathcal L_j w )(r) \Phi_j(\omega),\qquad\mathcal L_j =\partial_{rr} +\frac{d-1}{r}\partial_r -\frac{\lambda_j}{r^2}.
\end{equation}
The solutions of $\mathcal L_j w=0$ are linear combinations of $r^{\gamma_j^+}$ and $r^{-\gamma_j^-}$, where $\gamma_j^\pm\geq 0$ are given by
\begin{align}\label{eq:gammajpm}
\gamma_j^\pm &= \sqrt{\left(\frac{d-2}{2}\right)^2 +\lambda_j} \pm \frac{2-d}{2} \,,
\end{align}
that is,
\begin{align*}
\gamma_j^+ &= k \hspace{6.5em}  \text{for }\lambda_j  =k^2+k(d-2)\,,\\
\gamma_j^- &   = k +d-2 
\hspace{3em}  \text{for }\lambda_j 
 =k^2+k(d-2)\,.
\end{align*}
The decay rate $\gamma \geq d-2$ is fixed and we denote by $j_0=j_0(\gamma)$ the integer $j_0\geq 0$ such that
\begin{align*}
&\left\lbrace j\in \mathbb N \colon \gamma_j^- <\gamma\right\rbrace
=\lbrace 0,\ldots, j_0\rbrace,\\
& \left\lbrace j\in \mathbb N \colon \gamma_j^- \geq \gamma\right\rbrace
=\lbrace j_0 + 1, j_0 +2,\ldots \rbrace.
\end{align*}
%We observe that
%\begin{equation} \label{e.gamj0}
%\gamma_{j_0+1}^\pm =\gamma_{j_0}^{\pm} + 1.
%\end{equation}
%Any function $f\in L^2(\R^d)$ admits a spherical harmonics expansion
%\begin{equation*}
%f =\sum_{j\geq 0} f_j(r)\Phi_j(\omega),
%%\qquad \sum_{j\geq 0}\int_0^\infty f_j(r)^2 r^{d-1}\, dr <\infty,
%\end{equation*}
%and we set
%\begin{equation*}
%\Pi^+f =\sum_{j >j_0} f_j(r)\Phi_j(\omega),\qquad \Pi^-f = \sum_{0\leq j\leq j_0} f_j(r)\Phi_j(\omega).
%\end{equation*}
 The function  $f\in L^2(\R^d\setminus B_1)$ admits a spherical harmonics expansion
\begin{equation*}
f =\sum_{j\geq 0} f_j(r)\Phi_j(\omega),
\end{equation*}
and the decay assumption  on $f$ implies
\begin{align}\label{eq:fjdecay}
\sum_{j\geq 0}\int_R^\infty  f_j(r)^2r^{d+1}\, dr 
&
\lesssim \sum_{k\geq 0} (2^{k}R)^2\int_{2^k R \leq |x|\leq 2^{k+1}R}f^2 \, dx
\nonumber
\\
&
\lesssim 
\sum_{k\geq 0} 
(2^k R)^{d-2\gamma -2}\ln^{2\theta}(2^{k+1} R)
\nonumber
\\
&
\lesssim
 R^{d-2\gamma-2}
\ln^{2\theta}(2R)\,,
\end{align}
for all $R\geq 1$.
We define $u$ as 
\begin{align*}
u:=\sum_{j\geq 0} u_j(r)\Phi_j(\omega),
\end{align*}
where $u_j\in W^{2,2}_{loc}(0,\infty)$ satisfy
\begin{align*}
\mathcal L_j u_j =f_j.
\end{align*}
To write down an explicit formula for $u_j$ we rewrite $\mathcal L_j$, defined in \eqref{e.Lj}, as
\begin{align*}
\mathcal L_j u  =r^{-d+1+\gamma_j^-}\partial_r  [  r^{d-1-2\gamma_j^-}
\partial_r ( r^{\gamma_j^-}u) ],
\end{align*}
and define
\begin{align} \label{e.ujdef}
u_j(r)&=
\left\lbrace
\begin{aligned}
r^{-\gamma_j^-}\int_r^\infty t^{2\gamma_j^-+1-d}\int_t^\infty s^{d-1-\gamma_j^-} f_j(s)\, ds\, dt &\qquad\text{if }j\in\lbrace 0,\ldots, j_0\rbrace,\\
r^{-\gamma_j^-}\int_1^r t^{2\gamma_j^-+1-d}\int_t^\infty s^{d-1-\gamma_j^-} f_j(s)\, ds \, dt &\qquad\text{if }j \geq j_0 +1.
\end{aligned}
\right.
\end{align}
This is well defined because,
  for any $t \geq 1$,
   using Cauchy-Schwarz,  \eqref{eq:fjdecay}  with the choice $R = t$, and the fact that $\gamma_j^-\geq d-2>0$, 
we can estimate the inner integral by
\begin{align}\label{eq:estimfj1}
\int_t^\infty s^{d-1-\gamma_j^-} \abs{f_j(s)}\, ds
&\leq \left(\int_t^\infty s^{-2-2\gamma_j^-}s^{d-1}ds\right)^{\frac 12}
\left(\int_t^\infty s^2f_j(s)^2s^{d-1}ds\right)^{\frac 12}  
\nonumber
\\
&\lesssim
 \frac{1}{\sqrt{2\gamma_j^- + 2-d}}t^{\tfrac{d}{2} - \gamma_j^- - 1} t^{\tfrac{d}{2} - \gamma - 1}\ln^\theta(2t)
 \nonumber
 \\
 &
 = \frac{1}{\sqrt{ 2\gamma_j^- + 2-d}} t^{d - \gamma - \gamma_j^- -2}\ln^\theta(2t) . 
\end{align}
 Furthermore, as  $t\mapsto t^{2\gamma_j^- + 1 - d }t^{d-2-\gamma - \gamma_j^-}\ln^\theta t = t^{\gamma_j^--\gamma-1}\ln^\theta t$ is integrable near $\infty$ if $\gamma_j^-<\gamma$, i.e., if  $j\leq j_0$, the functions $u_j$ in \eqref{e.ujdef} are well-defined. 

Let $j\leq j_0$ and set
\begin{align*}
\alpha:=\gamma +\gamma_{j_0}^- +1 -d,
\end{align*}
so that  
$2\gamma +1 - d >\alpha > 2\gamma_j^-+1 - d$. By \eqref{eq:estimfj1} and Cauchy-Schwarz we have
\begin{align*}
|u_j(r)|^2 &
\leq \frac{r^{-2\gamma_j^-}}{2+2\gamma_j^--d}\left( 
\int_r^\infty t^{\gamma_j^--\frac{d}{2}} 
\left(\int_t^\infty s^2f_j(s)^2s^{d-1}ds\right)^{\frac 12}\, dt
 \right)^2 \\
& = \frac{r^{-2\gamma_j^-}}{2+2\gamma_j^--d}\left( 
\int_r^\infty t^{\gamma_j^--\frac{d}{2}-\frac{\alpha}{2}} 
t^{\frac{\alpha}{2}}\left(\int_t^\infty s^2f_j(s)^2s^{d-1}ds\right)^{\frac 12}\, dt
 \right)^2 \\
 & \leq 
 \frac{r^{-2\gamma_j^-}}{2+2\gamma_j^--d}
 \int_r^\infty t^{2\gamma_j^--d-\alpha}\, dt \int_r^{\infty}t^{\alpha} \left(\int_t^\infty s^2f_j(s)^2s^{d-1}ds\right)\, dt \\
 & = \frac{r^{-d+1-\alpha}}{(2+2\gamma_j^--d)(\alpha-2\gamma_j^- +d - 1)}
 \int_r^{\infty}t^{\alpha} \left(\int_t^\infty s^2f_j(s)^2s^{d-1}ds\right)\, dt\\
 &\leq 
  \frac{r^{-d+1-\alpha}}{(d-2)(\gamma-\gamma_{j_0}^-)}
   \int_r^{\infty}t^{\alpha} \left(\int_t^\infty s^2f_j(s)^2s^{d-1}ds\right)\, dt,
\end{align*}
where in the last line, we used that $\gamma_j^- \geqslant d-2$ so that $ 2 + 2\gamma_j^- - d \geqslant d-2,$ and that $\gamma + \gamma_{j_0}^- - 2\gamma_{j}^- \geqslant \gamma - \gamma_{j_0}^-,$ when $j \leqslant j_0.$
Summing and using \eqref{eq:fjdecay}, we deduce
\begin{align*}
\sum_{j=0}^{j_0} \frac{|u_j(r)|^2}{r^2} &
\leq \frac{r^{-d-1-\alpha}}{(d-2)(\gamma-\gamma_{j_0}^-)}
\int_r^{\infty}t^{\alpha} \left(\sum_{j=0}^{j_0}\int_t^\infty s^2f_j(s)^2s^{d-1}ds\right)\, dt \\
& \le  \frac{r^{-d-1-\alpha}}{(d-2)(\gamma-\gamma_{j_0}^-)} \int_r^{\infty}t^{\alpha +d-2\gamma -2}\ln^{2\theta}(2t)\, dt \\
& \lesssim
\frac{r^{-2\gamma-2}\ln^{2\theta}(2r)}{(d-2)(\gamma-\gamma_{j_0}^-)(2\gamma +1 -d -\alpha)}\\
& \leq \frac{r^{-2\gamma-2}\ln^{2\theta}(2r)}{(d-2)(\gamma-\gamma_{j_0}^-)^2}.
\end{align*}
For $j\geq j_0+1$
we need to distinguish cases if $\gamma=\gamma_{j_0+1}^-$.
We introduce $j_1  \geq j_0$ such that
\begin{align*}
\gamma 
&
=\gamma_j^-\quad\text{for }j\in\lbrace j_0+1,\ldots,j_1\rbrace,\\
\gamma
&
< \gamma_j^- \quad\text{for }j\geq j_1+1.
\end{align*}
For $j\geq j_1+1$
 we set 
\begin{align*}
\beta = \gamma +\gamma_{j_0+1}^- +1 -d,
\end{align*}
which satisfies $ 2\gamma +1 - d < \beta < 2\gamma_j^-+1 - d$.
Using \eqref{eq:estimfj1} and Cauchy-Schwarz we find
\begin{align*}
|u_j(r)|^2
&
\leq
\frac{r^{-2\gamma_j^-}}
{2+2\gamma_j^- -d}
\int_1^r t^{2\gamma_j^- -d - \beta}\, dt
 \int_1^{r}t^{\beta} \left(\int_t^\infty s^2f_j(s)^2s^{d-1}ds\right)\, dt
 \\
 &
 \lesssim
 \frac{
r^{-d +1 -\beta}
}
{(d-2)(\gamma-\gamma_{j_1+1}^-)}
 \int_1^{r}t^{\beta} \left(\int_t^\infty s^2f_j(s)^2s^{d-1}ds\right)\, dt
\end{align*}
so that
from \eqref{eq:fjdecay} 
we obtain that
\begin{align*}
\sum_{j=j_1+1}^{\infty} \frac{|u_j(r)|^2}{r^2}
& \lesssim \frac{r^{-2\gamma-2}\ln^{2\theta}(2r)}{(d-2)(\gamma_{j_1+1}^- -\gamma)^2}.
\end{align*}
It remains to treat $j_0+1\leq j\leq j_1$, where $\gamma=\gamma_j^-$.
In that case, the same manipulations, with $\beta=2\gamma+1-d$,
 lead to
\begin{align*}
|u_j(r)|^2
&
\leq
\frac{r^{-2\gamma}}
{d-2}
\ln r 
 \int_1^{r}t^{2\gamma +1 - d} \left(\int_t^\infty s^2f_j(s)^2s^{d-1}ds\right)\, dt
\,,
\end{align*}
and, using \eqref{eq:fjdecay},
\begin{align*}
\sum_{j=j_0+1}^{j_1} \frac{|u_j(r)|^2}{r^2}
& \lesssim \frac{r^{-2\gamma-2}\ln^{2\theta+2}(2r)}{d-2}.
\end{align*}
We conclude that
\begin{align*}
\sum_{j=0}^\infty
\frac{|u_j(r)|^2}{r^2}
& \lesssim \frac{1}{d-2}\left(\frac{1}{(\gamma-\gamma_{j_0}^-)^2}+\frac{1}{(\gamma_{j_1+1}^- -\gamma)^2} +\ln^2(2 r) \right) r^{-2\gamma -2}\ln^{2\theta}(2r)
\end{align*}
Therefore, since $\gamma\geq d-2$,
\begin{align*}
&
\frac{1}{R^d}\int_{|x|\geq R}\frac{|u|^2}{|x|^2}\, dx
=\frac{1}{R^d}\int_R^\infty\left(  \sum_{j=0}^\infty \frac{|u_j(r)|^2}{r^2}\right) \, r^{d-1}\, dr \\
& 
\lesssim \left(\frac{1}{(\gamma-\gamma_{j_0}^-)^2}+\frac{1}{(\gamma_{j_0+1}^- -\gamma)^2} +\ln^2(2R)\right) 
\frac{R^{-2\gamma-2}\ln^{2\theta}(2R)}{(d-2)^2}
\,,
\end{align*}
which implies the conclusion.
\end{proof}

%Lemma~\ref{l:decaylogfullspace} implies directly the following result:

%\begin{cor}\label{c:decaylog}
%Let $d\geq 3$, $\gamma\geq d-2$ and $\theta\geq 0$,  $R_*>\lambda \geq 1$ and  $f$ a function in $B_{R_*}\setminus B_\lambda$ satisfying 
%\begin{align*} 
%\left(\frac{1}{R^d}\int_{R<|x|<2R} f^2 \, dx  \right)^{\frac 12}&\leq \frac{\ln^\theta (2R)}{R^{\gamma+2}}\qquad\forall R\in [\lambda ,R_*/2].
%\end{align*}
%Then there exists a function $u$ such that $\Delta u =f$ in $B_{R_*}\setminus B_\lambda$ and
%\begin{align*} 
%\left(\frac{1}{R^d}\int_{R<|x|<2R}u^2\, dx\right)^{\frac 12}
% &\leq C \frac{\ln^{1+\theta}(2 R)}{R^{\gamma}}\qquad\forall R\in [\lambda,R_*/2],
%\end{align*}
%where $C>0$ depends only on $d, \gamma$ and $\theta$.
%\end{cor}

\section{Decay of harmonic functions in annuli}\label{a:harm}

In this appendix 
we gather
 two results
about 
pointwise control of
 harmonic functions in annuli.
 We use the same notations as in \S~\ref{s:harm_ext} and \S~\ref{a:poisson}, denoting by $\lbrace \Phi_j\rbrace$ an orthonormal system of eigenfunctions of the Laplacian on $\mathbb S^{d-1}$, with eigenvalues $\lambda_j$ and associated 
 powers $\gamma_j^\pm$ as in \eqref{eq:gammajpm}.
 
 \medskip
 
The first result is an annulus version of the fact that if a harmonic function in $\R^3\setminus B_\lambda$ is square-integrable, then it decays like $1/|x|^2$.

\begin{lem}\label{l:harm_L2_annulus}
Let $R_*/8 > \lambda \geq 1$.
Assume $\Delta u=0$ in $B_{R_*}\setminus B_\lambda \subset\R^3$ and
\begin{align*}
\int_{B_{R_*}\setminus B_\lambda}  |u|^2\, dx  \leq 1.
\end{align*}
Then we have
\begin{align*}
|u| +|x| |\nabla u|
\lesssim \frac{\sqrt\lambda}{|x|^2} 
+\frac{R_*^{-1/2}}{|x|} , 
\qquad\text{for }2\lambda\leq |x|\leq \frac{R_*}{4}\,.
\end{align*}
\end{lem}

The second result is an annulus version of the fact that a harmonic function with finite energy in $\R^d\setminus B_\lambda$ only has decaying modes.

\begin{lem}\label{l:harm_H1_annulus}
Let $R_*/8 > \lambda \geq 1$ and
assume $\Delta u=0$ in $B_{R_*}\setminus B_\lambda \subset\R^d$.
% and
%\begin{align*}
%\int_{B_{R_*}\setminus B_\lambda}|\nabla u|^2\, dx \leq 1\,.
%\end{align*}
Then we can decompose $u$ as a sum of two harmonic functions $u=v+w$, with
\begin{align*}
&
v(r\omega)=a_0\Phi_0 +\sum_{j\geq 0} b_j r^{-\gamma_j^-}\Phi_j(\omega)\,,
\\
&
\sum_{j\geq 0}\frac{|b_j|^2}{(4\lambda)^{2\gamma_j^-}}
\lesssim \frac{1}{\lambda^{d-2}} \int_{B_{2\lambda}\setminus B_\lambda}|\nabla u|^2\, dx\,,
\\
\text{and }
&
|w|^2+|x|^2|\nabla w|^2 \lesssim 
\frac{1}{R_*^{d-2}}
\int_{B_{R_*}\setminus B_{R_*/2}}
\!\!
|\nabla u|^2\, dx
\quad\text{for }2\lambda\leq |x|\leq \frac{R_*}{4}\,.
\end{align*}
The multiplicative constants depend on $d$.
\end{lem}

In the proofs of both lemmas, 
the main tool is an elementary estimate on the coefficients of 
a harmonic function generated by one single spherical harmonic, that is,
$u(r\omega) =\big( ar^{\gamma_k^+} +b r^{-\gamma_k^-} \big)\Phi_k(\omega)$,
 in terms of integrals of $|u|^2$.

\begin{lem}\label{l:control_ab}
For any $d\geq 2$, any $a,b\in\R$ and $\gamma^\pm \geq 0$, 
 such that
 $\gamma^+ +\gamma^-\geq 1$ and $\gamma^--\gamma^+ =d-2$, the function
\begin{align*}
u(r)=a \, r^{\gamma^+} +b r^{-\gamma^-}\qquad\text{for }r>0\,,
\end{align*}
satisfies, for any $\mu>1$, 
the estimates
\begin{align*}
|a|^2R^{2\gamma^+} &
\!\!
\lesssim 
\frac{\mu^{-2\gamma^+}}{R^d}
\int_{R}^{\mu^3 R} |u|^2\, r^{d-1}\, dr\,,
\\
|b|^2 R^{-2\gamma^-}
&
\!\!
\lesssim 
\frac{\mu^{4\gamma^-}}{R^d}\int_{R}^{\mu^3 R}|u|^2\, r^{d-1}\, dr
\,,\qquad\forall R>0\,,
\end{align*}
where the multiplicative constant depends on $\mu$ and $d$, but not on $\gamma^\pm$.
\end{lem}

\begin{proof}[Proof of Lemma~\ref{l:control_ab}]
We denote by $A$ the average of $|u|^2$ on $[R,\mu^3 R]$ with respect to $r^{d-1}\, dr$, that is,
\begin{align*}
A=\frac{d}{(\mu^{3d} -1) R^d}\int_{R}^{\lambda R}|u|^2\, r^{d-1}\, dr,\,.
\end{align*}
By the mean value theorem, we can find $R_1\in [R,\mu R]$ and $R_2\in [\mu^2 R,\mu^3 R]$ such that
\begin{align*}
|u(R_1)|^2
&
\leq
\frac{d}{(\mu^d -1)R^d} \int_{R}^{\mu R}|u|^2\, r^{d-1}\, dr =\frac{\mu^{3d} -1}{\mu^d -1} A\,,
\\
|u(R_2)|^2
&
\leq \frac{d}{\mu^{2d}(\mu^d-1)R^d} \int_{\mu^2 R}^{\mu^3 R}|u|^2\, r^{d-1}\, dr
=\frac{\mu^{3d} - 1}{\mu^{2d}(\mu^d-1)} A\,.
\end{align*}
Inverting the system
\begin{align*}
a\, R_1^{\gamma^+}\! +b \, R_1^{-\gamma^-} \!
&
=u(R_1)\,,\\
a\, R_2^{\gamma^+}\! +b \, R_2^{-\gamma^-} \!
&
=u(R_2)\,,
\end{align*}
we obtain
\begin{align*}
a
&
=
\frac{ -R_2^{-\gamma^-} u(R_1) 
+
R_1^{-\gamma^-}u(R_2) }
{R_2^{\gamma^+} R_1^{-\gamma^-} 
-
R_1^{\gamma^+}R_2^{-\gamma^-}
}
\, ,
\\
b
&
=
\frac{R_2^{\gamma^+}u(R_1) -  R_1^{\gamma^+} u(R_2) }
{R_2^{\gamma^+} R_1^{-\gamma^-} 
-
R_1^{\gamma^+}R_2^{-\gamma^-}
}
 \,. 
\end{align*}
Using that $R_1\leq \mu R\leq \mu^2 R\leq R_2$ 
 we find 
\begin{align*}
R_2^{\gamma^+} R_1^{-\gamma^-} -
R_1^{\gamma^+}R_2^{-\gamma^-}
\geq  
 \frac{\mu^{\gamma^+}-\mu^{-\gamma^-}}
{(\mu R)^{\gamma^--\gamma^+}}\,,
\end{align*}
and, using also that $R_1\in [R,\mu R]$, $R_2\in [\mu^2 R,\mu^3R]$,  
we deduce
\begin{align*}
|a|
&
\leq 
\frac{1+\mu^{-2\gamma^-}}{1-\mu^{-(\gamma^+ +\gamma^-)}}
\mu^{\gamma^- -\gamma^+}
\!
(\mu R)^{-\gamma^+}
\!\!
\max\big(|u(R_1)|,|u(R_2)|\big)
\\
|b|
&
\leq
\frac{1+\mu^{-2\gamma^+}}
{1-\mu^{-(\gamma^+ +\gamma^-)}}
\mu^{\gamma^+ -\gamma^-}
\!
(\mu^2 R)^{\gamma^-}
\!\!
\max\big(|u(R_1)|,|u(R_2)|\big)\,.
\end{align*}
Squaring these inequalities
and using that $|u(R_j)|^2\lesssim A$,
$\gamma^++\gamma^-\geq 1$
and  $\gamma^--\gamma^+ =d-2$,
we conclude.
\end{proof}

With Lemma~\ref{l:control_ab} now proven,
we establish the estimates of Lemma~\ref{l:harm_L2_annulus} and Lemma~\ref{l:harm_H1_annulus}.

\begin{proof}[Proof of Lemma~\ref{l:harm_L2_annulus}]
%By scaling, we assume without loss of generality $\lambda=1$. 
%(Otherwise consider the rescaled function $\tilde{u}(x) = \lambda^\frac32 u(x/\lambda)$.)
%
We write the spherical harmonics decomposition
\begin{align*}
u(r\omega) =\sum_{j\geq 0} u_j(r)\Phi_j(\omega),
\qquad u_j(r)=a_j r^{\gamma_j^+} +b_j r^{-\gamma_j^-}\,,
\end{align*}
and denote
\begin{align*}
A_j(R)=\frac{1}{R^3}\int_{R}^{2R} |u_j|^2\, r^2\, dr\,,
\qquad
\tilde A_j(R)=R^3 A_j(R)\,,
\end{align*}
so that
\begin{align}\label{eq:sumtildeAjL2}
\sum_{j\geq 0} \tilde A_j(R) \lesssim \int_{B_{2R}\setminus B_R}|u|^2\, dx \lesssim 1\qquad\forall R\in [\lambda,R_*/2]\,.
\end{align}
Applying Lemma~\ref{l:control_ab} to each $u_j$ with $\mu=2^{1/3}$ we have the inequalities
\begin{align*}
|a_j|^2
&
 \lesssim (2^{1/3}R)^{-2\gamma_j^+-3} \tilde A_j(R)\,,
\\
|b_j|^2
&
 \lesssim (2^{2/3}R)^{2\gamma_j^- - 3} \tilde A_j(R)\,,\qquad\forall R\in [\lambda,R_*/2]\,.
\end{align*}
Recall $\gamma_0^-=1$ and $\gamma_j^-\geq 2$ for $j\geq 1$, so the sign of the exponent of $R$ in the inequality for $b_j$ is different for $j=0$ and $j\geq 1$.
Choosing $R=R_*/2$ in the estimate on $a_j$ and $b_0$,
and $R=\lambda$ 
in the estimate on $b_j$ for $j\geq 1$, 
we obtain
\begin{align*}
|a_j|^2
&
 \lesssim 
 \frac{2^{\frac 4 3 \gamma_j^+  }}{R_*^{2\gamma_j^+ +3} } \tilde A_j(R_*/2)
 \qquad\text{for }j\geq 0\,,
 \\
 |b_0|^2
 &\lesssim 
 \frac{1}{R_*}\,,
\qquad
|b_j|^2
 \lesssim   
2^{\frac 43 \gamma_j^-}\lambda^{2\gamma_j^- -3} 
 \tilde A_j(\lambda)
 \qquad\text{for }j\geq 1\,.
\end{align*}
We use this  to obtain
\begin{align*}
\Xint{-}_{B_{2R}\setminus B_R}|u|^2\, dx
&
\lesssim \sum_{j\geq 0} \frac{1}{R^3}\int_{R}^{2R}|a_j r^{\gamma_j^+}+b_j r^{-\gamma_j^-}|^2\, r^2\, dr
\\
&
\lesssim \sum_{j\geq 0}|a_j|^2 (2R)^{2\gamma_j^+}
+\sum_{j\geq 0} |b_j|^2 R^{-2\gamma_j^-}
\\
&
\lesssim
\frac{1}{R_*^3} \sum_{j\geq 0}
\Big(\frac{2^{5/3}R}{R_*}\Big)^{2\gamma_j^+}
\tilde A_j(R_*/2)
+\frac{1/R_*}{R^2} 
\\
&\quad
+\frac{\lambda}{R^4}\sum_{j\geq 1}\Big(\frac{2^{2/3}\lambda}{R}\Big)^{2\gamma_j^- -4}\tilde A_j(\lambda),
\end{align*}
and therefore, using the summability property \eqref{eq:sumtildeAjL2} of the $\tilde A_j$'s,
\begin{align*}
\Xint{-}_{B_{2R}\setminus B_R}|u|^2\, dx
\lesssim \frac{1/R_*}{R^2} +\frac{\lambda}{R^4}
\qquad\text{for }2^{2/3}\lambda \leq R\leq \frac{R_*}{2^{5/3}}\,.
\end{align*}
Using elliptic estimates for the harmonic function $u$, this implies the 
conclusion of Lemma~\ref{l:harm_L2_annulus}.
\end{proof}

\begin{proof}[Proof of Lemma~\ref{l:harm_H1_annulus}]
We write the spherical harmonics decomposition
\begin{align*}
u(r\omega) =\sum_{j\geq 0} u_j(r)\Phi_j(\omega),
\qquad u_j(r)=a_j r^{\gamma_j^+} +b_j r^{-\gamma_j^-}\,,
\end{align*}
and denote
\begin{align*}
A_j(R)=
\frac{1}{R^{d}}
\int_{R}^{2R}
|u_j|^2\, r^{d-1}\, dr\,,
\qquad
\hat A_j(R)=R^{d-2} A_j(R)\,,
\end{align*}
so that
\begin{align}\label{eq:sumhatAjH1}
\sum_{j\geq 1} \hat A_j(R) 
&
\lesssim \int_{B_{2R}\setminus B_R}\frac{|\nabla_\omega u|^2}{|x|^2}\, dx 
\nonumber
\\
&
\lesssim \int_{B_{2R}\setminus B_R} |\nabla u|^2 \, dx 
\qquad\forall R\in [\lambda,R_*/2]\,.
\end{align}
Applying Lemma~\ref{l:control_ab} to each $u_j$ with $\mu=2^{1/3}$ we have the inequalities
\begin{align*}
|a_j|^2
&
 \lesssim (2^{1/3}R)^{-2\gamma_j^+-d+2} \hat A_j(R)\,,
\\
|b_j|^2
&
 \lesssim (2^{2/3}R)^{2\gamma_j^- - d+2} \hat A_j(R)\,,\qquad\forall R\in [\lambda,R_*/2]\,.
\end{align*}
Choosing $R=R_*/2$ in the estimate on $a_j$,
and $R=\lambda$ 
in the estimate on $b_j$, 
we obtain, for all $j\geq 1$,
\begin{align*}
|a_j|^2
&
 \lesssim 
 \Big(\frac{2^{\frac 23}}{R_*}\Big)^{2\gamma_j^+ +d -2}
 \hat A_j(R_*/2)
 \,,
 \\
|b_j|^2
&
 \lesssim   
\big(2^{\frac 23} \lambda\big)^{2\gamma_j^--d+2}
 \hat A_j(\lambda)
 \,.
\end{align*}
Using \eqref{eq:sumhatAjH1} we deduce
\begin{align*}
\sum_{j\geq 1}\frac{|b_j|^2}{(4\lambda)^{2\gamma_j^-}}
\lesssim  \frac{1}{\lambda^{d-2}}\sum_{j\geq 1}\hat A_j(\lambda)
\lesssim \frac{1}{\lambda^{d-2}}\int_{B_{2\lambda}\setminus B_\lambda}|\nabla u|^2\, dx\,.
\end{align*}
For $j=0$ we can use  that $\partial_r u_0 =-\gamma_0^- b_0 r^{-\gamma_0^- -1}$, 
to obtain
\begin{align*}
\frac{|b_0|^2}{(4\lambda)^{2\gamma_0^-}}
\lesssim \frac{1}{\lambda^{d-2}}\int_{B_{2\lambda}\setminus B_\lambda}|\partial_r u_0|^2\, dx
\lesssim \frac{1}{\lambda^{d-2}}\int_{B_{2\lambda}\setminus B_\lambda}|\nabla u|^2\, dx\,.
\end{align*}
This shows that the function $v$ given by
\begin{align*}
v(r\omega)&
=a_0\Phi_0 + \sum_{j\geq 0} b_j r^{-\gamma_j^-}\Phi_j(\omega)\,,
\end{align*}
does satisfy the claimed estimate.
It remains to prove the estimate on the function $w=u-v$ given by
\begin{align*}
w(r\omega) 
&
=\sum_{j\geq 1} a_j r^{\gamma_j^+}\Phi_j(\omega)\,.
\end{align*}
For any $R\in [\lambda,R_*/2^{2/3}]$, we use the above estimate on $a_j$ and the control \eqref{eq:sumhatAjH1} on the sum of the $\hat A_j$'s
to
 calculate
\begin{align*}
\Xint{-}_{B_{2R}\setminus B_R}|w|^2\, dx
&
\lesssim 
\sum_{j\geq 1} |a_j|^2 R^{2\gamma_j^+}
\\
&
\lesssim 
\frac{1}{R_*^{d-2}}\sum_{j\geq 1}\Big(\frac{2^{2/3}R}{R_*}\Big)^{2\gamma_j^+}\hat A_j(R_*/2)
\\
&
\lesssim \frac{1}{R_*^{d-2}}\sum_{j\geq 1}\hat A_j(R_*/2)
\\
&
\lesssim\frac{1}{R_*^{d-2}}\int_{B_{R_*}\setminus B_{R_*/2}}
\!\!
|\nabla u|^2\, dx\,.
\end{align*}
The conclusion follows from elliptic estimates for the harmonic function $w$.
\end{proof}

\end{appendices}
\bibliographystyle{acm}
\bibliography{interactions}

\end{document}